\documentclass{article}
\usepackage{amssymb}
\usepackage{amsfonts}
\usepackage{amsmath}

\newtheorem{theorem}{Theorem}

\newtheorem{corollary}[theorem]{Corollary}

\newtheorem{lemma}[theorem]{Lemma}

\newtheorem{remark}[theorem]{Remark}

\begin{document}

\title{On the prospective minimum of the random walk conditioned to stay
nonnegative\thanks{%
This work was supported by the Russian Science Foundation under grant
no.24-11-00037 https://rscf.ru/en/project/24-11-00037/ }}
\author{V.A.Vatutin\thanks{%
Steklov Mathematical Institute Gubkin street 8 119991 Moscow Russia Email:
vatutin@mi.ras.ru}, E.E.Dyakonova\thanks{%
Steklov Mathematical Institute Gubkin street 8 119991 Moscow Russia Email:
elena@mi-ras.ru}}
\maketitle

\begin{abstract}
Let
\begin{equation*}
S_{0}=0,\quad S_{n}=X_{1}+...+X_{n},\ n\geq 1,
\end{equation*}%
be a random walk whose increments belong without centering to the domain of
attraction of a stable law with scaling constants $a_{n}$, that provide
convergence as $n\rightarrow \infty $ of the distributions of the sequence $%
\left\{ S_{n}/a_{n},n=1,2,...\right\} $ to this stable law. Let $%
L_{r,n}=\min_{r\leq m\leq n}S_{m}$ \ be the minimum of the random walk on
the interval $[r,n]$. It is shown that
\begin{equation*}
\lim_{r,k,n\rightarrow \infty }\mathbf{P}\left( L_{r,n}\leq ya_{k}|S_{n}\leq
ta_{k},L_{0,n}\geq 0\right) ,\, t\in \left( 0,\infty \right),
\end{equation*}%
can have five different expressions, the forms of which depend on the
relationships between the parameters $r,k$ and $n$.

\textbf{Key words}: random walks, stable distributions, conditional limit
theorems
\end{abstract}

\section{Introduction}

We consider a random walk
\begin{equation*}
S_{0}=0,\quad S_{n}=X_{1}+...+X_{n},\ n\geq 1,
\end{equation*}%
whose increments belong to the domain of attraction of a stable law. To give
a more detailed description of the random walk class we are interesting in
introduce the set
\begin{equation*}
\mathcal{A}:=\{\alpha \in (0,2)\backslash \{1\},\,|\beta |<1\}\cup \{\alpha
=1,\beta =0\}\cup \{\alpha =2,\beta =0\}\subset \mathbb{R}^{2}.
\end{equation*}%
For a pair $(\alpha ,\beta )\in \mathcal{A}$ and a random variable $X$ we
write $X\in \mathcal{D}\left( \alpha ,\beta \right) $ if the distribution of
$X$ belongs to the domain of attraction of a stable law with density $%
g_{\alpha ,\beta }(x),x\in (-\infty ,+\infty ),$ and the characteristic
function%
\begin{equation*}
G_{\alpha ,\beta }(w)=\int_{-\infty }^{+\infty }e^{iwx}g_{\alpha .\beta
}(x)\,dx=\exp \left\{ -c|w|^{\,\alpha }\left( 1-i\beta \frac{w}{|w|}\tan
\frac{\pi \alpha }{2}\right) \right\} ,\ c>0,
\end{equation*}%
and, in addition, $\mathbf{E}X=0$ if this moment exists. This implies, in
particular, that there is an increasing sequence of positive numbers
\begin{equation}
a_{n}\ =\ n^{1/\alpha }\ell (n)  \label{DefA}
\end{equation}%
with a slowly varying at infinity sequence $\ell (1),\ell (2),\ldots ,$ such
that, as $n\rightarrow \infty $
\begin{equation}
\mathcal{L}\left\{ \frac{S_{\left[ nt\right] }}{a_{n}},t\geq 0\right\}
\Longrightarrow \mathcal{L}\left\{ \mathcal{Y}=\left\{ Y_{t},t\geq 0\right\}
\right\} ,  \label{dva}
\end{equation}%
where%
\begin{equation*}
\mathbf{E}e^{iwY_{t}}=G_{\alpha ,\beta }(wt^{1/\alpha }),t\geq 0,
\end{equation*}%
and the symbol $\Longrightarrow $ stands for the weak convergence in the
space $D[0,\infty )$ of c\`{a}dl\`{a}g functions endowed with Skorokhod
topology. Observe that if $X_{n}\overset{d}{=}X\in \mathcal{D}\left( \alpha
,\beta \right) $ for all $n\in \mathbb{N}:\mathbf{=}\left\{ 1,2,...\right\} $
then%
\begin{equation*}
\lim_{n\rightarrow \infty }\mathbf{P}\left( S_{n}>0\right) =\rho =\mathbf{P}%
\left( Y_{1}>0\right) \in (0,1).
\end{equation*}%
We assume throughout the paper that the random walk under consideration
meets one of the following restrictions.

\textbf{Condition A1.} \textit{The random variables }$X_{n},n\in \mathbb{N},$%
\textit{\ are independent copies of a random variable }$X\in \mathcal{D}%
\left( \alpha ,\beta \right) $\textit{.\ Besides, the distribution of }$X$%
\textit{\ is non-lattice.}\emph{\ }

\textbf{Condition A2. }\textit{The law of }$X$\textit{\ \ is absolutely
continuous with respect to the Lebesgue measure on }$\mathbb{R}$\textit{,
and there exists }$n\in \mathbb{N}$\textit{\ such that the density }$%
g_{n}(x):=\mathbf{P}(S_{n}\in dx)/dx$\textit{\ of }$S_{n}$\textit{\ is
bounded (therefore, }$g_{n}(x)\in L^{\infty }$\textit{).}

\begin{remark}
We remind that \ the requirement $g_{n}(x)\in L^{\infty }$ for some $n\in
\mathbb{N}$ is the standard necessary and sufficient condition for the
uniform convergence of the rescaled density $g_{n}(x/a_{n})/a_{n}$ to the
density of $Y_{1}$ (see, for instance, \cite[Section 46]{GK54}).
\end{remark}

\begin{remark}
Using the results for lattice distributions obtained in \cite{CC2013} and
following (with evident changes) the line of arguments we use below to
establish our main results one can prove the respective statements for the
lattice distributions satisfying
\end{remark}

\textbf{Condition A1'}. The law of $X$ belongs to $\mathcal{D}\left( \alpha
,\beta \right) $ and is supported on the lattice $\{0,\pm 1,\pm 2,...\}$
with maximal span 1.

We set%
\begin{eqnarray*}
L_{r,n} &=&\min_{r\leq m\leq n}S_{m},\quad L_{n}=L_{0,n},\quad
M_{r,n}=\max_{r\leq m\leq n}S_{m},\quad M_{n}=M_{1,n} \\
\tau _{r,n} &=&\min \left\{ r\leq m\leq n:S_{m}=L_{r,n}\right\} ,\quad \tau
_{n}=\tau _{0,n}=\min \left\{ 0\leq m\leq n:S_{m}=L_{n}\right\} .
\end{eqnarray*}

The aim of the paper is to investigate the asymptotic behavior of the
distribution of the random variable $L_{r,n}$ given $\left\{ S_{n}\leq
ta_{k},L_{n}\geq 0\right\} $ when $k=o(n)$ and $r=r(n)$ vary in an
appropriate way as $n\rightarrow \infty $ .

In what follows we denote by $\mathbf{P}_{w}\left( \cdot \right) $ the
distributions generated by the random walks stating from $S_{0}=w$ and use
the natural agreement $\mathbf{P}=\mathbf{P}_{0}\left( \cdot \right) .$

Random walks conditioned to stay positive, nonnegative or negative on some
set were investigated by many authors. Thus, Bolthausen \cite{Bolt76}
(improving the results of Iglehart \cite{Ig74}) has shown that if \ $\mathbf{%
E}X=0,\mathbf{E}X^{2}\in (0,\infty )$ then
\begin{equation}
\mathcal{L}\left\{ \frac{S_{\left[ nt\right] }}{\sigma \sqrt{n}},t\in
\lbrack 0,1]|L_{n}\geq 0\right\} \Longrightarrow \mathcal{L}\left\{ \mathcal{%
B}^{+}\right\} ,  \label{Meander2}
\end{equation}%
where $\mathcal{B}^{+}=\left\{ B_{t}^{+},t\in \lbrack 0,1]\right\} $ is a
Brownian meander. Durrett \cite{Dur78} extended this result by showing that
if $X\in \mathcal{D}\left( \alpha ,\beta \right) ,$ then
\begin{equation*}
\mathcal{L}\left\{ \frac{S_{\left[ nt\right] }}{a_{n}},t\in \lbrack
0,1]|L_{n}\geq 0\right\} \Longrightarrow \mathcal{L}\left\{ \mathcal{Y}%
^{+}\right\}
\end{equation*}%
where $\mathcal{Y}^{+}=\left\{ Y_{t}^{+},t\in \lbrack 0,1]\right\} $ is a
Levy meander. These results have been generalized in \cite{BD94}, \cite%
{CC2008} and \cite{CD2010} to the case when the random walk is conditioned
to stay positive on the set $[0,\infty )$.

Local versions of Bolthausen's and Durrett's theorems were obtained by
Caravenna \cite{Car2005} for the case $X\in \mathcal{D}\left( 2,0\right) $
and by Vatutin and Wachtel \cite{VW09} for the general case $X\in \mathcal{D}%
\left( \alpha ,\beta \right) $, respectively.

Relatively recently Caravenna and Chaumont \cite{CC2013} have made the next
step in studying such random walks. They impose different conditions on $%
S_{0}$ and $S_{n}$ and assume that the random walk remains nonnegative on $%
[0,n ]$. More precisely, given $n\in \mathbb{N}$ and $x,y\in \lbrack
0,\infty )$ they call a random walk starting at $S_{0}=x\geq 0$, ending at $%
S_{n}=y\geq 0$ and staying nonnegative (or positive) on the interval $\left[
1,n-1\right] $ a bridge of length $n$ and analysed the following laws%
\begin{equation}
\mathbf{P}_{x,y}^{\uparrow ,n}(\cdot )=\mathbf{P}_{x}(\cdot |L_{1,n-1}\geq
0,S_{n}=y),  \label{DefBridge1}
\end{equation}%
\begin{equation}
\mathbf{\hat{P}}_{x,y}^{\uparrow ,n}(\cdot )=\mathbf{P}_{x}(\cdot
|L_{1,n-1}>0,S_{n}=y).  \label{DefBridge2}
\end{equation}%
In order for the conditioning in the right-hand sides of (\ref{DefBridge1})
and (\ref{DefBridge2}) to be well-defined, they worked in the lattice case
under \ a condition slightly more general than Condition A1' and in the
absolutely continuous case under Conditions A1, A2 and the assumption $%
g_{n}^{+}(x,y)>0$, where

\begin{eqnarray*}
&&g_{n}^{+}(x,y):=\frac{\mathbf{P}_{x}\left( L_{1,n-1}>0,S_{n}\in dy\right)
}{dy}  \notag \\
&=&\int_{\mathcal{K}(n-1)}\left[ g(s_{1}-x)\left(
\prod_{i=2}^{n-1}g(s_{i}-s_{i-1})\right) g(y-s_{n-1})\right]
ds_{1}...ds_{n-1}  \label{ContDistrib}
\end{eqnarray*}%
with $\mathcal{K}(n-1):=\left\{ s_{1}>0,...,s_{n-1}>0\right\} $ and $g(\cdot
)=g_{\alpha ,\beta }(\cdot )$ the density of the distribution of the
increments of the random walk. It was shown in \cite{CC2013} that if
\begin{equation}
x/a_{n}\rightarrow 0,y/a_{n}\rightarrow 0\text{ as }n\rightarrow \infty
\label{XYneglig}
\end{equation}%
then%
\begin{equation}
\mathcal{L}\left\{ \frac{S_{\left[ nt\right] }}{a_{n}},t\in \lbrack
0,1]|S_{0}=x,L_{1,n-1}\geq 0,S_{n}=y\right\} \Longrightarrow \mathcal{L}%
\left\{ \mathcal{Y}^{++}\right\} ,  \label{TotalCVonvergence}
\end{equation}%
where $\mathcal{Y}^{++}=\left\{ Y_{t}^{++},t\in \lbrack 0,1]\right\} $ is a
Levy bridge (excursion) to stay positive with $\mathbf{P}\left(
Y_{0}^{++}=Y_{1}^{++}=0\right) =1$ (see Section 6 in \cite{CC2013} and \cite%
{Urb2011} for more detail concerning the definition of the process).

If $X\in \mathcal{D}\left( 2,0\right) $ we get in the limit a standard
Brownian excursion (see, for instance, \cite{ItoMaC74}, p.75). It is not
difficult to deduce from (\ref{TotalCVonvergence}) and theorem 5.1 in \cite%
{CC2013}\ that given (\ref{XYneglig}) the following relations are valid as $%
n\rightarrow \infty $%
\begin{equation}
\mathcal{L}\left\{ \frac{S_{\left[ nt\right] }}{a_{n}},t\in \lbrack
0,1]|S_{0}=x,L_{1,n}\geq 0,S_{n}\leq y\right\} \Longrightarrow \mathcal{L}%
\left\{ \mathcal{Y}^{++}\right\}  \label{TotalConveregence2}
\end{equation}%
and
\begin{equation*}
\lim_{n\rightarrow \infty }\mathbf{P}\left( \frac{S_{n-m}}{a_{n}}\geq
z|S_{0}=x,L_{1,n-1}\geq 0,S_{n}=y\right) =0
\end{equation*}%
for any $z>0$ if $m=o(n)$. Therefore, $S_{n-m}/a_{n}\rightarrow 0$ in
probability as $n\rightarrow \infty $ and (\ref{TotalCVonvergence}) provides
practically no information about the behavior of $S_{n-m}$ in this case. A
more detailed description of the distribution of the random variable $%
S_{n-m} $ was given in \cite{VDD2023}, where it was shown that, as $%
n\rightarrow \infty $
\begin{equation}
\mathcal{L}\left\{ \frac{S_{n-m}-S_{n}}{a_{m}}|S_{0}=x,L_{1,n}\geq
0,S_{n}\leq y\right\} \rightarrow \mathcal{L(J)}  \label{SmallDeviat}
\end{equation}%
under the condition (\ref{XYneglig}), where the forms of the distribution of
the random variable $\mathcal{J}$ are different for the cases when $%
\lim_{m\rightarrow \infty }y/a_{m}$ is zero, is infinity or is a positive
constant. These results were used in \cite{VDD2023} to investigate the
asymptotic behavior of the population size distribution of a critical
branching process evolving in non-favorable random environment.

The next natural step in studying properties of random walks is to
investigate the behavior of various functionals specified on the tragectory $%
\{S_{r},S_{r+1},...,S_{n}\}$. For instance, if $\mathbf{E}X=0,\mathbf{E}%
X^{2}\in (0,\infty )$ then, in view of (\ref{Meander2}) \ and the
Donsker-Prokhorov invariance principle for continuous functionals
\begin{equation*}
\mathbf{P}\left( \frac{L_{\left[ nt\right] ,n}}{\sigma \sqrt{n}}\leq
x|L_{n}\geq 0\right) \rightarrow \mathbf{P}\left( \min_{t\leq s\leq
1}B_{s}^{+}\leq x\right)
\end{equation*}%
for any $t\in \lbrack 0,1]$ and $x>0$. This limiting relation plays an
important role in studying properties of the so-called reduced braniching
processes in random environment (see \cite{BV97}, \cite{Vat2002}).
Unfortunately, neither the functional limit theorem (\ref{TotalConveregence2}%
) nor the weak convergence (\ref{SmallDeviat}) can be used to directly
obtain meaningful results of such a kind for the properties of the important
for applications random variable $L_{r,n}=\min_{r\leq k\leq n}S_{k}$ under
the assumption $\left\{ S_{0}=x,L_{1,n}\geq 0,S_{n}\leq y\right\} $ if $%
y=o(a_{n})$ as $n\rightarrow \infty $.

We fill this gap in the present paper and investigate, for fixed $t\in
\left( 0,\infty \right) $ the asymptotic behavior of the distribution%
\begin{equation}
\mathbf{P}\left( L_{r,n}\leq x|S_{n}\leq ta_{k},L_{n}\geq 0\right)
\label{Distrib_minimum}
\end{equation}%
for $k=k(n)=o(n),r=r(n)\in \lbrack 0,n]$ and $\min \left\{ r,k\right\}
\rightarrow \infty $ as $n\rightarrow \infty $ using conditional local limit
theorems for random walks obtained in \cite{CC2013}, \cite{Don12} and \cite%
{VW09}. It happens that, depending on the relationship between $k,r$ and $n$
the limit of the distribution (\ref{Distrib_minimum}) has 5 different forms.

We consider these cases separately. For the convenience of readers we give
the list of sections devoted to the respective cases:

\bigskip

\frame{.$%
\begin{array}{ccc}
\text{Section \ref{Sec2}} & \text{-- the case } & n\gg k\gg r; \\
\text{Section \ref{Sec3}} & \text{-- the case } & n\gg k=\theta r,\quad
\theta \in (0,\infty ); \\
\text{Section \ref{Sec4}} & \text{-- the case } & \min (r,n-r)\gg k; \\
\text{Section \ref{Sec5}} & \text{-- the case } & n\gg k=\theta (n-r),\theta
\in (0,\infty ); \\
\text{Section \ref{Sec6}} & \text{-- the case } & n\gg k\gg n-r.%
\end{array}%
$}

\bigskip

Before to pass to the analyse of these cases we collect in Section \ref{Sec1}
a number of basic results which will be used in the proofs.

\section{Some properties of random walks\label{Sec1}}

In the sequel we denote by $C_{1},C_{2},...$ some absolute constants which
may not be the same in different formulas or even within the same
complicated formula.

We set $\mathbb{N}_{0}:=\mathbb{N}\mathbf{\cup }\left\{ 0\right\} $. Given
two positive sequences $\left\{ c_{n},n\in \mathbb{N}\right\} $, $\left\{
d_{n},n\in \mathbb{N}\right\} ,$ we write as usual $c_{n}\sim d_{n}$ if $%
lim_{n\rightarrow \infty }c_{n}/d_{n}=1;c_{n}=o(d_{n})$ or $c_{n}\ll d_{n}$
if $lim_{n\rightarrow \infty }c_{n}/d_{n}=0$; $c_{n}=O(d_{n})$ if $\lim
\sup_{n\rightarrow \infty }c_{n}/d_{n}<\infty $.

We recall that a sequence $\left\{ c_{n},n\in \mathbb{N}\right\} $ of
positive numbers or a real positive function $c(x)$ --- is said to be
regularly varying at infinity with index $\gamma \in \mathbb{R}$ , denoted $%
c_{n}\in R_{\gamma }$ or $c(x)\in R_{\gamma }$ if $c_{n}\sim n^{\gamma }l(n)$
$(c(x)\sim x^{\gamma }l(x))$, where $l(x)$ is a slowly varying function,
i.e. a positive real function with the property that $l(cx)/l(x)\rightarrow
1 $ as $x\rightarrow \infty $ for all fixed $c>0$.

Set $S_{0}:=0,$ $\tau _{0}^{\pm }:=0,$ and for $k\geq 1~$denote by%
\begin{equation*}
\tau _{k}^{-}:=\inf \left\{ n>\tau _{k-1}^{-}:S_{n}\leq S_{\tau
_{k-1}^{-}}\right\}
\end{equation*}%
the weak descending ladder variables and by%
\begin{equation*}
\tau _{k}^{+}:=\inf \left\{ n>\tau _{k-1}^{+}:S_{n}\geq S_{\tau
_{k-1}^{+}}\right\}
\end{equation*}
the weak ascending ladder variables of the sequence $S_{0},S_{1},...$. \ Put
\begin{equation*}
H_{k}^{\pm }:=\pm S_{\tau _{k}^{\pm }}
\end{equation*}%
and introduce the parameter%
\begin{equation*}
\zeta =\mathbf{P}\left( H_{1}^{+}=0\right) =\mathbf{P}\left(
H_{1}^{-}=0\right) \in (0,1),
\end{equation*}%
where to justify this equality it is necessary to use the duality principle
for random walks (see, for instance, \cite[Ch.XII, Sec.2]{Fel}):
\begin{equation*}
\left\{ S_{n}-S_{n-k},k=0,1,...,n\right\} \overset{d}{=}\left\{
S_{k},k=0,1,...,n\right\} .
\end{equation*}

For $x\geq 0$ introduce renewal functions%
\begin{equation*}
V^{\pm }(x)=\sum_{k=0}^{\infty }\mathbf{P}\left( H_{k}^{\pm }\leq x\right)
=\sum_{k=0}^{\infty }\sum_{n=0}^{\infty }\mathbf{P}\left( \tau _{k}^{\pm
}=n,\pm S_{n}\leq x\right) .\quad
\end{equation*}%
Note that
\begin{equation*}
V^{+}(x)=\sum_{n=0}^{\infty }\mathbf{P}\left( S_{n}\leq x,L_{n}\geq 0\right)
.
\end{equation*}%
It is easy to check that $V^{\pm }(x)$ are non-decreasing, right-continuous
and%
\begin{equation}
V^{\pm }(0)=\sum_{k=0}^{\infty }\mathbf{P}\left( H_{k}^{\pm }=0\right) =%
\frac{1}{1-\zeta }.  \label{V_zero}
\end{equation}%
It is known (see, for instance, \cite{Rog1971}, \cite{Sin57}) that
\begin{equation}
\mathbf{P}\left( \tau _{1}^{+}>n\right) =\mathbf{P}\left( M_{n}<0\right) \in
R_{-\rho },\quad V^{+}(x)\in R_{\alpha \rho },  \label{Regular1}
\end{equation}

\begin{equation}
\mathbf{P}\left( \tau _{1}^{-}>n\right) =\mathbf{P}\left( L_{1,n}>0\right)
\in R_{-(1-\rho )},\quad V^{-}(x)\in R_{\alpha (1-\rho )}.  \label{Regular2}
\end{equation}

Note that, in view of (\ref{V_zero})-(\ref{Regular2}) and the asymptotical
representation (31) in \cite{VW09} and (3.18) in \cite{CC2013} there are
positive constants $\hat{C}$, $\mathcal{C}^{+}$ and \ $\mathcal{C}^{-}$ such
that
\begin{equation*}
n\mathbf{P}\left( \tau _{1}^{-}>n\right) \mathbf{P}\left( \tau
_{1}^{+}>n\right) \sim \hat{C},
\end{equation*}%
\begin{equation}
V^{+}(a_{n})\sim \frac{\mathcal{C}^{+}}{1-\zeta }n\mathbf{P}\left( \tau
_{1}^{-}>n\right) ,\quad V^{-}(a_{n})\sim \frac{\mathcal{C}^{-}}{1-\zeta }n%
\mathbf{P}\left( \tau _{1}^{+}>n\right) ,  \label{AsymV_plus_minus}
\end{equation}%
and, therefore, as $n\rightarrow \infty $%
\begin{equation*}
\mathbf{P}\left( \tau _{1}^{+}>n\right) V^{+}(a_{n})\sim C^{\ast }:=\frac{%
\mathcal{C}^{+}\hat{C}}{1-\zeta }\in (0,\infty ),
\end{equation*}%
\begin{equation}
\mathbf{P}\left( \tau _{1}^{-}>n\right) V^{-}(a_{n})\sim C^{\ast \ast }:=%
\frac{\mathcal{C}^{-}\hat{C}}{1-\zeta }\in (0,\infty ).  \label{Product11}
\end{equation}

Let $g^{+}(\cdot )$ be the density of the time-one marginal distribution of
the meander of the Levy process $\mathcal{Y}$ (see \cite{Chau97}) and $%
g^{-}(\cdot )$ be the density of the time-one marginal distribution of the
meander of $-\mathcal{Y}$. As it is shown in \cite[Lemma 3]{VDD2023}
\begin{equation*}
\frac{1}{C^{\ast }}=\int_{0}^{\infty }z^{\alpha \rho }g^{-}(z)dz.
\end{equation*}%
Hence, by symmetry arguments we conclude that%
\begin{equation*}
\frac{1}{C^{\ast \ast }}=\int_{0}^{\infty }z^{\alpha (1-\rho )}g^{+}(z)dz.
\end{equation*}

For the construction of renewal functions the strict ladder variables $%
\left\{ \hat{\tau}_{k}^{\pm },k\geq 0\right\} $ and epochs $\left\{ \hat{H}%
_{k}^{\pm },k\geq 0\right\} $ are used along with the weak ladder variables.
They are defined as $\hat{\tau}_{0}^{\pm }:=0,$\ $\hat{H}_{0}^{\pm }:=0$
and, for $k\geq 1$
\begin{equation*}
\hat{\tau}_{k}^{\pm }:=\inf \left\{ n>\hat{\tau}_{k-1}^{\pm }:\pm S_{n}>\pm
S_{\hat{\tau}_{k-1}^{\pm }}\right\} ,\quad \hat{H}_{k}^{\pm }:=\pm S_{\tau
_{k}^{\pm }}.
\end{equation*}

The sequences $\left\{ \hat{H}_{k}^{\pm },k\geq 0\right\} $ and $\left\{
\hat{\tau}_{k}^{\pm },k\geq 0\right\} $ generate the renewal functions%
\begin{equation*}
\hat{V}^{\pm }(x):=\sum_{k=0}^{\infty }\mathbf{P}\left( \hat{H}_{k}^{\pm
}\leq x\right) =\sum_{k=0}^{\infty }\sum_{n=0}^{\infty }\mathbf{P}\left(
\hat{\tau}_{k}^{\pm }=n,\pm S_{n}\leq x\right) .
\end{equation*}%
It is easy to check that
\begin{equation*}
\hat{V}^{-}(x):=\sum_{n=0}^{\infty }\mathbf{P}\left( S_{n}\geq
x,M_{n}<0\right) .
\end{equation*}%
The connection between the introduced renewal functions is very simple (see
\cite[ch. XII.1 equation (1.13)]{Fel}):%
\begin{equation}
\hat{V}^{\pm }(x)=\left( 1-\zeta \right) V^{\pm }(x).  \label{RenewRelation}
\end{equation}

In the sequel to avoid complicated notation we suppose as a rule that the
distribution of the increments $X_{i},i=1,2,...,$ is absolutely continuous.
This means, in particulary, that $\hat{V}^{\pm }(x)=V^{\pm }(x)$ for all $%
x\geq 0.$\ This agreement allows us to use without further explanations the
results from other papers, that deal with the functions $\hat{V}^{\pm }(x),$
$V^{\pm }(x),$ as well as the functions
\begin{equation*}
\underline{{V}}^{\pm }(x):=\sum_{k=0}^{\infty }\mathbf{P}\left( H_{k}^{\pm
}<x\right) =\sum_{k=0}^{\infty }\sum_{n=0}^{\infty }\mathbf{P}\left( \tau
_{k}^{\pm }=n,\pm S_{n}<x\right)
\end{equation*}%
and
\begin{equation*}
\underline{\hat{V}}^{\pm }(x):=\sum_{k=0}^{\infty }\mathbf{P}\left( \hat{H}%
_{k}^{\pm }<x\right) =\sum_{k=0}^{\infty }\sum_{n=0}^{\infty }\mathbf{P}%
\left( \hat{\tau}_{k}^{\pm }=n,\pm S_{n}<x\right) .
\end{equation*}

\bigskip The fundamental property of $\hat{V}^{-}$is the identity%
\begin{equation*}
\mathbf{E}[\hat{V}^{-}(x+X_{1});X_{1}+x\geq 0]\ =\ \hat{V}^{-}(x)\
\end{equation*}%
which, in view of $\hat{V}^{-}(0)=1$ implies%
\begin{equation}
\mathbf{E}[\hat{V}^{-}(S_{n});L_{n}\geq 0]\ =1  \label{harm2}
\end{equation}%
for all $n\geq 1$.

In the proofs presented below we need the following asymptotic
represenations, following from (\ref{Regular1}), (\ref{Regular2}) and the
properties of regularly varying functions:

(1) as $x\rightarrow \infty ,$
\begin{eqnarray}
\left( \alpha \rho +1\right) \int_{0}^{x}V^{+}(w)dw &\sim &xV^{+}(x),  \notag
\\
&&  \label{AsympV} \\
\quad \left( \alpha (1-\rho )+1\right) \int_{0}^{x}V^{-}(w)dw &\sim
&xV^{-}(x);  \notag
\end{eqnarray}

(2) as $n\rightarrow \infty $
\begin{eqnarray}
\int_{b}^{c}V^{+}(za_{n})g^{-}(z)dz &\sim &V^{+}(a_{n})\int_{b}^{c}z^{\alpha
\rho }g^{-}(z)dz,  \label{IntegralV_plus} \\
\int_{b}^{c}V^{-}(za_{n})g^{+}(z)dz &\sim &V^{-}(a_{n})\int_{b}^{c}z^{\alpha
(1-\rho )}g^{+}(z)dz  \label{IntegralV_minus}
\end{eqnarray}
uniformly for $0\leq b<c$ from any fixed interval $[0,C].$

Denote%
\begin{equation}
b_{n}=\frac{1}{na_{n}}=\frac{1}{n^{1+1/\alpha }\ell (n)}.  \label{Defb}
\end{equation}

In the sequel we will use several times the following simple observation.

\begin{lemma}
\label{L_estimB}If $a_{n},n=1,2,...,$ satisfies (\ref{DefA}), then there
exists a constant $C\in (0,\infty )$ such that
\begin{equation}
\sum_{j=k}^{n-k}b_{j}b_{n-j}\leq C\frac{b_{n}}{a_{k}}  \label{B_sum}
\end{equation}%
for all $n>2k.$
\end{lemma}

\textbf{Proof}. Since the sequence $b_{m},m=1,2,...,$ is decreasing, it
follows that
\begin{equation*}
\sum_{j=k}^{n-k}b_{j}b_{n-j}\leq b_{\left[ n/2\right] }\left(
\sum_{j=[n/2]}^{n-k}b_{n-j}+\sum_{j=k}^{[n/2]}b_{j}\right) \leq 2b_{\left[
n/2\right] }\sum_{j=k}^{\infty }b_{j}.
\end{equation*}%
for $n>2k$. \ Properties of regularly varying functions and (\ref{Defb})
allow us to conclude that there is a constant $C_{1}$ such that%
\begin{equation*}
b_{\left[ n/2\right] }\leq C_{1}b_{n}
\end{equation*}%
for all $n=3,4,...$ and (see \cite[Ch. VIII, Sec.9, Theorem1]{Fel}) that, as
$k\rightarrow \infty $%
\begin{equation}
\sum_{j=k}^{\infty }b_{j}=\sum_{j=k}^{\infty }\frac{1}{j^{1+1/\alpha }\ell
(j)}\sim \alpha kb_{k}=\frac{\alpha }{a_{k}}.  \label{Tauber1}
\end{equation}%
Combining the obtained estimates proves Lemma \ref{L_estimB}.

The next lemma contains several inequalities showing importance of the
renewal functions we introduced above.

\begin{lemma}
\label{L_Renewal_negl}If $X_{1}\in \mathcal{D}(\alpha ,\beta )$ then

(1) for any $\varepsilon >0$ there exists $K_{0}=K_{0}(\varepsilon )$ such
that
\begin{equation}
\sum_{j=K}^{\infty }\mathbf{P}\left( S_{j}\leq x,L_{j}\geq 0\right) \leq
\varepsilon V^{+}(x)  \label{Renewal_negl}
\end{equation}%
for all $K\geq K_{0}$ and all $x\geq 0;$

(2) for any $\varepsilon >0$ there exists $K_{0}=K_{0}(\varepsilon )$ such
that
\begin{equation}
\sum_{j=K}^{\infty }\mathbf{P}\left( S_{j}\geq -x,M_{j}<0\right) \leq
\varepsilon V^{-}(x)  \label{Renewal_negl1}
\end{equation}%
for all $K\geq K_{0}$ and all $x\geq 0;$

(3) there exists a constant $C\in (0,\infty )$ such that%
\begin{equation}
\mathbf{P}_{w}\left( S_{n}\in \lbrack x,x+y),L_{n}\geq 0\right) \leq
Cb_{n}V^{-}(w)\int_{x}^{x+y}V^{+}(u)du  \label{IntrervalEstimate}
\end{equation}%
for all nonnegative $w,x$ and $y$.

(4) there exists a constant $C\in (0,\infty )$ such that%
\begin{equation}
\mathbf{P}_{w}\left( S_{n}\in \lbrack -x,-x+y),M_{n}<0\right) \leq
Cb_{n}V^{+}(-w)\int_{x-y}^{x}V^{-}(u)du  \label{IntrervalEstimate1}
\end{equation}%
for all negative $w$ and $x\geq y\geq 0$.
\end{lemma}

\textbf{Proof}. It follows from (\ref{RenewRelation}) and Proposition 2.3 in
\cite{agkv} (rewritten in the notation of the present paper) that there
exists a constant $C$ such that, for all $n\geq 1$ and $w,x\geq 0$
\begin{equation}
\mathbf{P}_{w}\left( x<S_{n}\leq x+1,L_{n}\geq 0\right) \leq
Cb_{n}V^{-}(w)V^{+}(x).  \label{EstLocal}
\end{equation}%
Therefore,%
\begin{equation*}
\sum_{j=K}^{\infty }\mathbf{P}\left( S_{j}\leq x,L_{j}\geq 0\right) \leq
CV^{+}(x)\sum_{j=K}^{\infty }b_{j}\leq \frac{C_{1}}{a_{K}}V^{+}(x)
\end{equation*}%
proving (\ref{Renewal_negl}).

The estimate (\ref{IntrervalEstimate}) follows from (\ref{EstLocal}) by
integration over the interval $[x,x+y)$.

To prove the inequalies\ (\ref{Renewal_negl1}) and (\ref{IntrervalEstimate1}%
) it is necessary to use the estimate
\begin{equation*}
\mathbf{P}_{w}\left( x<S_{n}\leq x+1,M_{n}<0\right) \leq
Cb_{n}V^{+}(-w)V^{-}(-x)  \label{EstLocal1}
\end{equation*}%
valid for some constant $C$ and all \ $n\geq 1$ and $w,x\leq 0$, and
following from (\ref{RenewRelation}) and Proposition 2.3 in \cite{agkv}.

The lemma is proved.

Let%
\begin{equation*}
\mathcal{B}(x,n):=\{S_{n}\leq x,L_{n}\geq 0\}.
\end{equation*}%
According to Corollary 2 in \cite{VD2022} (written in the notation of the
present paper)
\begin{equation}
\mathbf{P}_{w}\left( \mathcal{B}(x,n)\right) \sim g_{\alpha ,\beta
}(0)V^{-}(w)b_{n}\int_{0}^{x}\underline{{V}}^{+}(u)du
\label{BasicAsymptotic0}
\end{equation}%
as $n\rightarrow \infty $ uniformly in \thinspace $x,w\geq 0$ such that $%
\max (x,w)\leq \delta _{n}a_{n},$ where $\delta _{n}\rightarrow 0$ as $%
n\rightarrow \infty $. In the sequel we investigate the case $x\rightarrow
\infty .$ In this situation we may change (\ref{BasicAsymptotic0}) by the
asymptotic relation
\begin{equation}
\mathbf{P}_{w}\left( \mathcal{B}(x,n)\right) \sim g_{\alpha ,\beta
}(0)V^{-}(w)b_{n}\int_{0}^{x}V^{+}(u)du  \label{BasicAsymptotic}
\end{equation}%
valid as $n\rightarrow \infty $ uniformly in \thinspace $x,w\geq 0$ such
that $0\leq \min (x,w)\leq \max (x,w)\leq \delta _{n}a_{n},$ where $\delta
_{n}\rightarrow 0$ as $n\rightarrow \infty $.

\begin{lemma}
\label{L_Levy}If condition A1 is valid, $n>>k$ then, for any $0<c<1,$ $w\geq
0$ and $t>0$
\begin{equation*}
\lim_{\varepsilon \downarrow 0,N\uparrow \infty }\lim_{cn\geq r\rightarrow
\infty }\mathbf{P}_{w}\left( S_{r}\notin \lbrack \varepsilon a_{r},Na_{r}]|%
\mathcal{B}(ta_{k},n)\right) =0.
\end{equation*}
\end{lemma}

\textbf{Proof}. Using (\ref{IntrervalEstimate}),(\ref{RenewRelation}), the
estimate $b_{n-r}\leq Cb_{n}$ (valid for $r\leq C_{1}n$ in view of (\ref%
{Defb})) and (\ref{BasicAsymptotic})\ we obtain%
\begin{eqnarray*}
\mathbf{P}_{w}\left( S_{r}\notin \lbrack \varepsilon a_{r},Na_{r}],\mathcal{B%
}(ta_{k},n)\right) &=&\int_{s\notin \lbrack \varepsilon a_{r},Na_{r}]}%
\mathbf{P}_{w}\left( S_{r}\in ds,L_{r}\geq 0\right) \mathbf{P}_{s}\left(
\mathcal{B}(ta_{k},n-r)\right) \\
&\leq &Cb_{n-r}\int_{0}^{ta_{k}}V^{+}(z)dz\int_{s\notin \lbrack \varepsilon
a_{r},Na_{r}]}V^{-}(s)\mathbf{P}_{w}\left( S_{r}\in ds,L_{r}\geq 0\right) \\
&\leq &C_{1}b_{n}\int_{0}^{ta_{k}}V^{+}(z)dz\int_{s\notin \lbrack
\varepsilon a_{r},Na_{r}]}\hat{V}^{-}(s)\mathbf{P}_{w}\left( S_{r}\in
ds,L_{r}\geq 0\right) \\
&\leq &C_{2}\mathbf{P}_{w}\left( \mathcal{B}(ta_{k},n)\right) \int_{s\notin
\lbrack \varepsilon a_{r},Na_{r}]}\hat{V}^{-}(s)\mathbf{P}_{w}\left(
S_{r}\in ds,L_{r}\geq 0\right) .
\end{eqnarray*}%
Clearly,%
\begin{eqnarray*}
\int_{s\notin \lbrack \varepsilon a_{r},Na_{r}]}\hat{V}^{-}(s)\mathbf{P}%
_{w}\left( S_{r}\in ds,L_{r}\geq 0\right) &=&\mathbf{E}_{w}\left[ \hat{V}%
^{-}(S_{r})I\left\{ S_{r}\notin \lbrack \varepsilon a_{r},Na_{r}]\right\}
,L_{r}\geq 0\right] \\
&=&\hat{V}^{-}(w)\mathbf{\hat{P}}_{w}^{-}\left( S_{r}\notin \lbrack
\varepsilon a_{r},Na_{r}]\right) ,
\end{eqnarray*}%
where, for $w\geq 0$%
\begin{eqnarray*}
\mathbf{\hat{P}}_{w}^{-}(A) &=&\frac{1}{\hat{V}^{-}(w)}\mathbf{E}_{w}\left[
\hat{V}^{-}(S_{r})I\left\{ A\right\} ,L_{r}\geq 0\right] \\
&=&\mathbf{\hat{P}}_{w}^{-}(A)=\frac{1}{V^{-}(w)}\mathbf{E}_{w}\left[
V^{-}(S_{r})I\left\{ A\right\} ,L_{r}\geq 0\right]
\end{eqnarray*}%
is a probability measure in view of (\ref{harm2}).

It follows from \cite[Theorem 1.1]{CC2008} that, for any $w\geq 0$ and $y>0$%
\begin{equation}
\lim_{r\rightarrow \infty }\mathbf{\hat{P}}_{w}^{-}\left( S_{r}\leq
ya_{r}\right) =\mathbf{P}\left( Y_{1}^{+}\leq y\right) ,  \label{Meand}
\end{equation}%
where, as before, $\left\{ Y_{s}^{+},0\leq s\leq 1\right\} $ is a Levy
meander with index $\alpha $. Since the distribution of the random variable $%
Y_{1}^{+}$ has no atoms at zero and infinity, the statement of the lemma
follows now from (\ref{Meand}).

\begin{corollary}
\label{C_conditNegligible}If condition A1 is valid, $n\gg k$ then, for any
event $\mathcal{G}_{n}$ measurable with respect to the $\sigma -$algebra
generated by the sequence $S_{0},S_{1},...,S_{n}$ and any $w\geq 0,\,0<c<1$
and $t>0$
\begin{equation*}
\lim_{\varepsilon \downarrow 0,N\uparrow \infty }\lim_{cn\geq r\rightarrow
\infty }\mathbf{P}_{w}\left( \mathcal{G}_{n},S_{r}\notin \lbrack \varepsilon
a_{r},Na_{r}]|\mathcal{B}(ta_{k},n)\right) =0.
\end{equation*}
\end{corollary}

\textbf{Proof}. The statement of the corollary follows from the inequality
\begin{equation*}
\mathbf{P}_{w}\left( \mathcal{G}_{n},S_{r}\notin \lbrack \varepsilon
a_{r},Na_{r}]|\mathcal{B}(ta_{k},n)\right) \leq \mathbf{P}_{w}\left(
S_{r}\notin \lbrack \varepsilon a_{r},Na_{r}]|\mathcal{B}(ta_{k},n)\right)
\end{equation*}%
and Lemma \ref{L_Levy}.

We now prove a lemma giving an upper estimate for sums connected with
renewal functions.

\begin{lemma}
\label{L_truncated}Let Condition A1 be valid. Then there exists $%
K_{0}=K_{0}(\varepsilon )$ such that for all nonnegative $Z>Y\geq 0$ and all
$K\geq K_{0}$
\begin{equation*}
\sum_{j=K}^{\infty }\mathbf{P}\left( S_{j}\in \lbrack Y,Z),L_{j}\geq
0\right) \leq \frac{C}{a_{K}}\int_{Y}^{Z}V^{+}(u)du
\end{equation*}%
and%
\begin{equation}
\sum_{j=K}^{\infty }\mathbf{P}\left( S_{j}\in \lbrack -Z,-Y),M_{j}<0\right)
\leq \frac{C}{a_{K}}\int_{Y}^{Z}V^{-}(u)du.  \label{new1}
\end{equation}
\end{lemma}

\textbf{Proof}. We prove (\ref{new1}) only.\textbf{\ }In view of (\ref%
{IntrervalEstimate1})
\begin{equation*}
\mathbf{P}\left( S_{j}\in \lbrack -Z,-Y),M_{j}<0\right) \leq
Cb_{j}\int_{Y}^{Z}V^{-}(u)du.
\end{equation*}%
Using (\ref{Tauber1}) we have
\begin{eqnarray*}
\sum_{j=K}^{\infty }\mathbf{P}\left( S_{j}\in \lbrack -Z,-Y),M_{j}<0\right)
&\leq &C\int_{Y}^{Z}V^{-}(u)du\sum_{j=K}^{\infty }b_{j} \\
&\leq &\frac{C_{1}}{a_{K}}\int_{Y}^{Z}V^{-}(u)du,
\end{eqnarray*}%
as required.

\begin{corollary}
\label{C_RenewalEquival}Let Condition A1 be valid. If $Z=o(a_{K})$ as $%
K\rightarrow \infty $ then
\begin{equation}
\sum_{j=0}^{K}\mathbf{P}\left( S_{j}\leq Z,L_{j}\geq 0\right) \sim V^{+}(Z)
\label{Tail11}
\end{equation}%
and%
\begin{equation}
\sum_{j=0}^{K}\mathbf{P}\left( S_{j}\geq -Z,M_{j}<0\right) \sim V^{-}(Z).
\label{Tail12}
\end{equation}
\end{corollary}

\textbf{Proof}. We have
\begin{eqnarray*}
0 &\leq &V^{+}(Z)-\sum_{j=0}^{K}\mathbf{P}\left( S_{j}\leq Z,L_{j}\geq
0\right) =\sum_{j=K+1}^{\infty }\mathbf{P}\left( S_{j}\leq Z,L_{j}\geq
0\right) \\
&\leq &\frac{C}{a_{K}}\int_{0}^{Z}V^{+}(u)du\leq \frac{CZ}{a_{K}}%
V^{+}(Z)=o(V^{+}(Z))
\end{eqnarray*}%
proving (\ref{Tail11}). Estimate (\ref{Tail12}) can be established by the
same arguments.

To write the relations given below in the more compact form we introduce the
notation $a\wedge b=\min \left( a,b\right) ,a\vee b=\max \left( a,b\right) .$

\begin{lemma}
\label{L_minLeft}Let Condition A1 be valid. If $n\gg k\gg r\rightarrow
\infty ,$ then for any positive $z,y$ and $t$

\begin{eqnarray}
\mathbf{P}_{za_{r}}(S_{\tau _{n}} &\leq &ya_{r},\mathcal{B}(ta_{k},n))
\notag \\
&=&\left( 1+o(1)\right) \mathbf{P}(\mathcal{B}(ta_{k},n))\times \left(
V^{-}(za_{r})-V^{-}((z-y\wedge z)a_{r})+\Delta _{1}\right) ,  \notag \\
&&  \label{Small_r}
\end{eqnarray}%
where%
\begin{equation*}
0\leq \Delta _{1}\leq \frac{C}{a_{k}}\left( 1+\int_{(z-y\wedge
z)a_{r}}^{za_{r}}V^{-}(u)du\right) .
\end{equation*}
\end{lemma}

\textbf{Proof}. Using the duality principle for random walks we have a
decomposition
\begin{eqnarray*}
&&\mathbf{P}_{za_{r}}\left( 0\leq S_{\tau _{n}}\leq ya_{r},S_{n}\leq
ta_{k}\right) =\mathbf{P}_{za_{r}}\left( 0\leq S_{\tau _{n}}\leq \left(
y\wedge z\right) a_{r},S_{n}\leq ta_{k}\right) \\
&=&\sum_{j=0}^{n}\int_{0}^{\left( y\wedge z\right) a_{r}}\mathbf{P}%
_{za_{r}}\left( S_{j}\in ds,\tau _{j}=j\right) \mathbf{P}\left( S_{n-j}\leq
ta_{k}-s,L_{n-j}\geq 0\right) \\
&\sim &\sum_{j=0}^{n}\mathbf{P}\left( S_{j}\in \left[ -za_{r},\left( \left(
y\wedge z\right) -z\right) a_{r}\right] ,M_{j}<0\right) \mathbf{P}\left(
S_{n-j}\leq ta_{k},L_{n-j}\geq 0\right)
\end{eqnarray*}%
as $n\gg k\gg r\rightarrow \infty ,$ since in this case $ta_{k}-s\sim ta_{k}$
for $0\leq s\leq za_{r}$. Without loss of generality we can assume that
\thinspace $n>2k.$ This estimate and point (4) of Lemma \ref{L_Renewal_negl}
imply
\begin{eqnarray*}
&&\sum_{j=k}^{n}\mathbf{P}\left( S_{n-j}\leq ta_{k},L_{n-j}\geq 0\right)
\mathbf{P}\left( S_{j}\in \left[ -za_{r},-za_{r}+\left( y\wedge z\right)
a_{r}\right] ,M_{j}<0\right) \\
&\leq &C\int_{(z-\left( y\wedge z\right)
)a_{r}}^{za_{r}}V^{-}(u)du\sum_{j=k}^{n}b_{j}\mathbf{P}\left( S_{n-j}\leq
ta_{k},L_{n-j}\geq 0\right) \\
&\leq &Ca_{r}V^{-}(za_{r})\sum_{j=k}^{n}b_{j}\mathbf{P}\left( S_{n-j}\leq
ta_{k},L_{n-j}\geq 0\right) .
\end{eqnarray*}%
By Lemma \ref{L_estimB} and point (3) of Lemma \ref{L_Renewal_negl}
\begin{eqnarray*}
&&\sum_{j=k}^{n-k}b_{j}\mathbf{P}\left( S_{n-j}\leq ta_{k},L_{n-j}\geq
0\right) \leq C\int_{0}^{ta_{k}}V^{+}(u)du\sum_{j=k}^{n-k}b_{j}b_{n-j} \\
&&\qquad \qquad \qquad \qquad \leq C\frac{b_{n}}{a_{k}}%
\int_{0}^{ta_{k}}V^{+}(u)du\leq \frac{C_{2}}{a_{k}}\mathbf{P}\left( \mathcal{%
B}(ta_{k},n)\right) .
\end{eqnarray*}

Further, using (\ref{Regular2}) and (\ref{AsymV_plus_minus}) and the
condition $n>>k$ we see that
\begin{eqnarray*}
&&\sum_{j=n-k}^{n}b_{j}\mathbf{P}\left( S_{n-j}\leq ta_{k},L_{n-j}\geq
0\right) \leq Cb_{n}\sum_{j=n-k}^{n}\mathbf{P}\left( S_{n-j}\leq
ta_{k},L_{n-j}\geq 0\right) \\
&&\qquad \leq Cb_{n}\sum_{q=0}^{k}\mathbf{P}\left( L_{q}\geq 0\right) \leq
C_{1}b_{n}k\mathbf{P}\left( L_{k}\geq 0\right) \sim C_{2}b_{n}V^{+}(a_{k}) \\
&&\qquad \qquad \qquad \qquad \leq \frac{C_{3}b_{n}}{a_{k}}%
\int_{0}^{a_{k}}V^{+}(u)du\leq \frac{C_{4}}{a_{k}}\mathbf{P}\left( \mathcal{B%
}(ta_{k},n)\right) .
\end{eqnarray*}%
Thus,%
\begin{equation}
\mathbf{P}_{za_{r}}\left( 0\leq S_{\tau _{n}}\leq ya_{r},S_{n}\leq
ta_{k};\tau _{n}\in \lbrack k,n]\right) \leq \frac{C_{2}}{a_{k}}\mathbf{P}%
\left( \mathcal{B}(ta_{k},n)\right) .  \label{Remainder1}
\end{equation}%
Observe that in view of (\ref{new1})
\begin{eqnarray*}
&&\sum_{j=0}^{k}\mathbf{P}\left( S_{n-j}\leq ta_{k},L_{n-j}\geq 0\right)
\mathbf{P}\left( S_{j}\in \left[ -za_{r},\left( y\wedge z-z\right) a_{r}%
\right] ,M_{j}<0\right) \\
&&\qquad \sim \mathbf{P}\left( S_{n}\leq ta_{k},L_{n}\geq 0\right)
\sum_{j=0}^{k}\mathbf{P}\left( S_{j}\in \left[ -za_{r},(y\wedge z-z)a_{r}%
\right] ,M_{j}<0\right) \\
&&\qquad \qquad =\mathbf{P}\left( S_{n}\leq ta_{k},L_{n}\geq 0\right) \left(
\hat{V}^{-}(za_{r})-\hat{V}^{-}((z-y\wedge z)a_{r})\right) \\
&& \\
&&\qquad \qquad \qquad \qquad \qquad \qquad \qquad \qquad \qquad +\mathbf{P}%
\left( S_{n}\leq ta_{k},L_{n}\geq 0\right) \Delta ,
\end{eqnarray*}%
where%
\begin{equation*}
0\leq \Delta \leq \frac{C_{1}}{a_{k}}\int_{(z-\left( y\wedge z\right)
)a_{r}}^{za_{r}}V^{-}(u)du.
\end{equation*}%
.

Combining the obtained estimates with (\ref{Remainder1}) we get (\ref%
{Small_r}).

Lemma \ref{L_minLeft} is proved.

The preparatory work we have done allows us to proceed in Sections 3-7 to
the statements and proofs of theorems describing the asymptotic behavior of
the probability $\mathbf{P}_{w}\left( S_{\tau _{r,n}}\leq ya_{r}|\mathcal{B}%
(ta_{k},n)\right) $, $t>0$, under various assumptions concerning the growth
rate of the parameters $k=k(n)$ and $r=r(n)$ as $n\rightarrow \infty $.

\section{The case $n\gg k\gg r$\label{Sec2}}

In this section we study the distribution of $S_{\tau _{r,n}}$ under the
condition $\{S_{n}\leq x,L_{n}\geq 0\}$ and relatively small $r\rightarrow
\infty $.

\begin{theorem}
\label{T_beginningNew} If Conditions A1 and A2 are valid then, for any $y>0$
and $t>0$
\begin{equation}
\lim_{n\gg k\gg r\rightarrow \infty }\mathbf{P}_{w}\left( S_{\tau
_{r,n}}\leq ya_{r}|\mathcal{B}(ta_{k},n)\right) =C^{\ast \ast }\mathcal{H}%
(y).  \label{CloseBegin}
\end{equation}%
where $C^{\ast \ast }$ is the same as in (\ref{Product11}) and
\begin{equation*}
\mathcal{H}(y)=\int_{0}^{\infty }g^{+}\left( z\right) \left( z^{\alpha
(1-\rho )}-(z-\left( y\wedge z\right) )^{\alpha (1-\rho )}\right) dz.
\end{equation*}
\end{theorem}

\begin{remark}
Observe that in view of (\ref{Product11}) the right-hand side of (\ref%
{CloseBegin}) tends to 1 as $y\rightarrow \infty $ and, therefore, the
limiting distribution specified by the right-hand side of (\ref{CloseBegin})
is proper.
\end{remark}

\textbf{Proof of Theorem \ref{T_beginningNew}}. We start by the asymptotic
representation%
\begin{equation}
\mathbf{P}_{w}\left( S_{r}\in dz,L_{r}\geq 0\right) =\frac{\mathbf{P}\left(
\tau _{1}^{-}>r\right) }{a_{r}}V^{-}(w)\left( g^{+}\left( \frac{z}{a_{r}}%
\right) +o(1)\right) dz  \notag
\end{equation}%
valid for all $z\in \lbrack 0,\infty )$ (see \cite[formula (5.1)]{CC2013}).
Observing that $\inf_{\varepsilon \leq z\leq N}g^{+}(z)>0$ for fixed $%
N>\varepsilon >0,$ we write the chain of relations
\begin{eqnarray}
&&\mathbf{P}_{w}\left( S_{\tau _{r,n}}\leq ya_{r},\mathcal{B}%
(ta_{k},n)\right)  \notag \\
&=&\int_{0}^{\infty }\mathbf{P}_{w}\left( S_{r}\in dz,L_{r}\geq 0\right)
\mathbf{P}_{z}\left( 0\leq S_{\tau _{n-r}}\leq ya_{r},S_{n-r}\leq
ta_{k}\right) dz  \notag \\
&=&\frac{\mathbf{P}\left( \tau _{1}^{-}>r\right) }{a_{r}}V^{-}(w)\int_{0}^{%
\infty }\left( g^{+}\left( \frac{z}{a_{r}}\right) +o(1)\right) \mathbf{P}%
_{z}\left( 0\leq S_{\tau _{n-r}}\leq ya_{r},S_{n-r}\leq ta_{k}\right) dz
\notag \\
&=&\mathbf{P}\left( \tau _{1}^{-}>r\right) V^{-}(w)\int_{0}^{\infty }\left(
g^{+}\left( z\right) +o(1)\right) \mathbf{P}_{za_{r}}\left( 0\leq S_{\tau
_{n-r}}\leq \left( y\wedge z\right) a_{r},S_{n-r}\leq ta_{k}\right) dz
\notag \\
&\sim &\mathbf{P}\left( \tau _{1}^{-}>r\right) V^{-}(w)\int_{\varepsilon
}^{N}g^{+}\left( z\right) \mathbf{P}_{za_{r}}\left( 0\leq S_{\tau
_{n-r}}\leq \left( y\wedge z\right) a_{r},S_{n-r}\leq ta_{k}\right) dz
\notag \\
&&+\mathbf{P}_{w}\left( S_{\tau _{r,n}}\leq ya_{r},S_{r}\notin \lbrack
\varepsilon a_{r},Na_{r}],\mathcal{B}(ta_{k},n)\right) ,   \notag
\end{eqnarray}%
where, according to Corollary \ref{C_conditNegligible}
\begin{equation*}
\lim_{\varepsilon \downarrow 0,N\uparrow \infty }\lim_{n\gg k\gg
r\rightarrow \infty }\frac{\mathbf{P}_{w}\left( S_{\tau _{r,n}}\leq
ya_{r},S_{r}\notin \lbrack \varepsilon a_{r},Na_{r}],\mathcal{B}%
(ta_{k},n)\right) }{\mathbf{P}_{w}\left( \mathcal{B}(ta_{k},n)\right) }=0.
\end{equation*}%
Using Lemma \ref{L_minLeft} with $n$ replaced by $n-r$ we obtain%
\begin{eqnarray*}
&&\int_{\varepsilon }^{N}g^{+}\left( z\right) \mathbf{P}_{za_{r}}\left(
0\leq S_{\tau _{n-r}}\leq ya_{r},S_{n-r}\leq ta_{k}\right) dz \\
&\sim &\mathbf{P}(\mathcal{B}(ta_{k},n-r))\int_{\varepsilon }^{N}g^{+}\left(
z\right) \left( V^{-}(za_{r})-V^{-}((z-\left( y\wedge z\right)
)a_{r})+\Delta _{1}\right) dz.
\end{eqnarray*}%
Since $n\gg k\gg r$, relations (\ref{BasicAsymptotic}) and (\ref{Defb})
imply
\begin{equation*}
\mathbf{P}(\mathcal{B}(ta_{k},n-r))\sim \mathbf{P}(\mathcal{B}(ta_{k},n))
\end{equation*}%
as $r\rightarrow \infty $. Further, in view of (\ref{Regular2}) and (\ref%
{IntegralV_minus}) we have
\begin{eqnarray*}
0 &\leq &\int_{\varepsilon }^{N}g^{+}\left( z\right) \Delta _{1}dz\leq \frac{%
C}{a_{k}}\int_{\varepsilon }^{N}g^{+}\left( z\right) \left(
1+\int_{(z-\left( y\wedge z\right) )a_{r}}^{za_{r}}V^{-}(u)du\right) dz \\
&\leq &\frac{C_{1}a_{r}}{a_{k}}\int_{\varepsilon }^{N}g^{+}\left( z\right)
\int_{z-\left( y\wedge z\right) }^{z}V^{-}(qa_{r})dqdz \\
&\leq &C_{2}\frac{a_{r}V^{-}(a_{r})}{a_{k}}\int_{\varepsilon
}^{N}g^{+}\left( z\right) \int_{z-\left( y\wedge z\right) }^{z}q^{\alpha
(1-\rho )}dqdz \\
&\leq &C_{3}\frac{a_{r}V^{-}(a_{r})}{a_{k}}\int_{0}^{N}g^{+}\left( z\right)
z^{\alpha (1-\rho )+1}dz=C_{4}\frac{a_{r}V^{-}(a_{r})}{a_{k}},
\end{eqnarray*}%
where $C_{4}$ is a constant depending on $N$ only. Taking into account (\ref%
{Product11}) we conclude that if $k\gg r\rightarrow \infty $, then for fixed
$N>\varepsilon >0$%
\begin{eqnarray*}
&&\mathbf{P}\left( \tau _{1}^{-}>r\right) V^{-}(w)\int_{\varepsilon
}^{N}g^{+}\left( z\right) \Delta _{1}dz\leq C_{4}V^{-}(w)\mathbf{P}\left(
\tau _{1}^{-}>r\right) V^{-}(a_{r})\frac{a_{r}}{a_{k}} \\
&\leq &C_{5}C^{\ast \ast }V^{-}(w)\frac{a_{r}}{a_{k}}\rightarrow 0
\end{eqnarray*}%
as $k>>r\rightarrow \infty $. Note now that in view of (\ref{IntegralV_minus}%
) and (\ref{Product11})
\begin{eqnarray*}
&&\int_{\varepsilon }^{N}g^{+}\left( z\right) \left( \hat{V}^{-}(za_{r})-%
\hat{V}^{-}((z-y\wedge z)a_{r})\right) dz \\
&&\qquad =\hat{V}^{-}(a_{r})\int_{\varepsilon }^{N}g^{+}(z)\left( \frac{\hat{%
V}^{-}(za_{r})-\hat{V}^{-}((z-y\wedge z)a_{r})}{\hat{V}^{-}(a_{r})}\right) dz
\\
&&\qquad \sim \hat{V}^{-}(a_{r})\int_{\varepsilon }^{N}g^{+}(z)\left(
z^{\alpha (1-\rho )}-(z-y\wedge z)^{\alpha (1-\rho )}\right) dz
\end{eqnarray*}%
as $r\rightarrow \infty $. Using (\ref{Product11}) once again we see that%
\begin{eqnarray*}
&&\mathbf{P}\left( \tau _{1}^{-}>r\right) V^{-}(w)\int_{\varepsilon
}^{N}g^{+}\left( z\right) \mathbf{P}_{za_{r}}\left( 0\leq S_{\tau
_{n-r}}\leq ya_{r},S_{n-r}\leq ta_{k}\right) dz \\
&\sim &C^{\ast \ast }V^{-}(w)\mathbf{P}(\mathcal{B}(ta_{k},n))\int_{%
\varepsilon }^{N}g^{+}\left( z\right) \left( z^{\alpha (1-\rho )}-(z-\left(
y\wedge z\right) )^{\alpha (1-\rho )}\right) dz.
\end{eqnarray*}%
Combining the obtained estimates we conclude that
\begin{equation*}
\mathbf{P}_{w}\left( S_{\tau _{r,n}}\leq ya_{r}|\mathcal{B}(ta_{k},n)\right)
\rightarrow C^{\ast \ast }\int_{0}^{\infty }g^{+}\left( z\right) \left(
z^{\alpha (1-\rho )}-(z-\left( y\wedge z\right) )^{\alpha (1-\rho )}\right)
dz
\end{equation*}%
as \ $n\gg k\gg r\rightarrow \infty $.

Theorem \ref{T_beginningNew} is proved.

\section{The case $n\gg k=\protect\theta r$\label{Sec3}}

In this section we assume that \ $n\gg k=\left[ \theta r\right] \rightarrow
\infty $ and, to avoid cumbersome formulas agree to consider $\sqrt{nr}$ as $%
\left[ \sqrt{nr}\right] $ and $\theta r$ as $\left[ \theta r\right] .$

\begin{theorem}
\label{T_initialDimilar} Let Conditions A1 and A2 be valid and $k=\left[
\theta r\right] ,\theta \in (0,\infty )$. Then, for any fixed $w\geq 0,t>0$
any $y\in \left[ 0,t\theta ^{1/\alpha }\right] $
\begin{equation}
\lim_{n\gg k\rightarrow \infty }\mathbf{P}_{w}\left( S_{\tau _{r,n}}\leq
ya_{r}|\mathcal{B}(ta_{k},n)\right) =W(t\theta ^{1/\alpha },y),
\label{Propor1}
\end{equation}%
where%
\begin{eqnarray*}
W(t,y) &=&\frac{C^{\ast \ast }(\alpha \rho +1)}{t^{\alpha \rho +1}}%
\int_{0}^{\infty }g^{+}\left( z\right) dz\int_{(z-y)\vee 0}^{z}q^{\alpha
\rho }(t-z+q)^{\alpha (1-\rho )}dq \\
&&+\frac{C^{\ast \ast }\alpha (1-\rho )}{t^{\alpha \rho +1}}\int_{0}^{\infty
}g^{+}\left( z\right) dz\int_{(z-y)\vee 0}^{z}\left( t-z+q\right) ^{\alpha
\rho +1}q^{\alpha (1-\rho )-1}dq.
\end{eqnarray*}
\end{theorem}

\textbf{Proof}. We first consider the asymptotic behavior of $\mathbf{P}%
_{w}\left( S_{\tau _{r,n}}\leq ya_{r},\mathcal{B}(ta_{r},n)\right) $ and
then, using the equivalence%
\begin{equation}
a_{k}=a_{\theta r}\sim \theta ^{1/\alpha }a_{r}  \label{Proport0}
\end{equation}%
and the inclusions%
\begin{equation}
\mathcal{B}(t\theta ^{1/\alpha }(1-\delta )a_{r},n)\subset \mathcal{B}%
(ta_{k},n)\subset \mathcal{B}(t\theta ^{1/\alpha }(1+\delta )a_{r},n),
\label{Inclusion}
\end{equation}%
valid for any $\delta \in (0,1)$ prove (\ref{Propor1}).

The same as in Theorem \ref{T_beginningNew} for fix $0<\varepsilon <N<\infty
$ we have
\begin{eqnarray*}
&&\mathbf{P}_{w}\left( S_{\tau _{r,n}}\leq ya_{r},\mathcal{B}%
(ta_{r},n)\right) \\
&\sim &\mathbf{P}\left( \tau _{1}^{-}>r\right) V^{-}(w)\int_{\varepsilon
}^{N}g^{+}\left( z\right) \mathbf{P}_{za_{r}}\left( 0\leq S_{\tau
_{n-r}}\leq \left( y\wedge z\right) a_{r},S_{n-r}\leq ta_{r}\right) dz \\
&&+\mathbf{P}_{w}\left( S_{\tau _{r,n}}\leq ya_{r},S_{r}\notin \lbrack
\varepsilon a_{r},Na_{r}],\mathcal{B}(ta_{r},n)\right)
\end{eqnarray*}%
as $n>>k\rightarrow \infty ,$ where
\begin{equation*}
\lim_{\varepsilon \downarrow 0,N\uparrow \infty }\lim_{n\gg r\rightarrow
\infty }\frac{\mathbf{P}_{w}\left( S_{\tau _{r,n}}\leq ya_{r},S_{r}\notin
\lbrack \varepsilon a_{r},Na_{r}],\mathcal{B}(ta_{r},n)\right) }{\mathbf{P}%
_{w}\left( \mathcal{B}(ta_{r},n)\right) }=0
\end{equation*}%
according to Corollary \ref{C_conditNegligible}. We consider the
decomposition
\begin{eqnarray*}
&&\mathbf{P}_{za_{r}}\left( 0\leq S_{\tau _{n-r}}\leq \left( y\wedge
z\right) a_{r},S_{n-r}\leq ta_{r}\right) \\
&=&\mathbf{P}\left( -za_{r}\leq S_{\tau _{n-r}}\leq \left( \left( y-z\right)
\wedge 0\right) a_{r},S_{n-r}\leq (t-z)a_{r}\right) \\
&=&Q_{z}\left( 0,\sqrt{nr}\right) +Q_{z}\left( \sqrt{nr}+1,n-\sqrt{nr}%
\right) +Q_{z}\left( n-\sqrt{nr}+1,n-r\right) ,
\end{eqnarray*}%
where, for $A<B$%
\begin{eqnarray*}
Q_{z}(A,B) &=&\mathbf{P}_{za_{r}}\left( 0\leq S_{\tau _{n-r}}\leq \left(
y\wedge z\right) a_{r},S_{n-r}\leq ta_{r},\tau _{n-r}\in \lbrack A,B]\right)
\\
&=&\sum_{j=A}^{B}\int_{-za_{r}}^{((y-z)\wedge 0)a_{r}}\mathbf{P}\left(
S_{j}\in ds,\tau _{j}=j\right) \mathbf{P}\left( S_{n-r-j}\leq
(t-z)a_{r}-s,L_{n-r-j}\geq 0\right) \\
&=&\sum_{j=n-k-B}^{n-k-A}\int_{-za_{r}}^{((y-z)\wedge 0)a_{r}}\mathbf{P}%
\left( S_{n-r-j}\in ds,\tau _{n-r-j}=n-r-j\right) \\
&&\times \mathbf{P}\left( S_{j}\leq (t-z)a_{r}-s,L_{j}\geq 0\right) .
\end{eqnarray*}%
First observe that by the duality principle for random walks the estimates
\begin{eqnarray*}
&&\mathbf{P}\left( S_{n-r-j}\in ds,\tau _{n-r-j}=n-r-j\right) \\
&=&\mathbf{P}\left( S_{n-r-j}\in ds,\tau _{1}^{+}>n-r-j\right) \sim
g_{\alpha ,-\beta }(0)b_{n-j-r}V^{-}(-s)ds \\
&\sim &g_{\alpha ,-\beta }(0)b_{n}a_{r}V^{-}(-qa_{r})dq
\end{eqnarray*}%
are valid for $r\ll n,$ $j\leq \sqrt{nr}$ uniformly in $q\in \lbrack
-z,(y-z)\wedge 0)$ (compare with Theorem 5.1 in \cite{CC2013}, where the
respective statement was proved for a random walk conditioned to stay
nonnegative). Since $a_{r}=o\left( a_{\sqrt{nr}}\right) $ as $n\rightarrow
\infty $ we have by (\ref{Tail11}) that
\begin{eqnarray*}
&&Q_{z}(n-\sqrt{nr}+1,n-r) \\
&=&\sum_{j=0}^{\sqrt{nr}}\int_{-z}^{(y-z)\wedge 0}\mathbf{P}\left(
S_{n-r-j}\in a_{r}dq,\tau _{n-r-j}=n-r-j\right) \\
&&\qquad \qquad \times \mathbf{P}\left( S_{j}\leq \left( t-z-q\right)
a_{r},L_{j}\geq 0\right) \\
&&\,\sim g_{\alpha ,\beta }(0)b_{n}a_{r}\int_{-z}^{(y-z)\wedge
0}V^{-}(-qa_{r})\left( \sum_{j=0}^{\sqrt{nr}}\mathbf{P}\left( S_{j}\leq
\left( t-z-q\right) a_{r},L_{j}\geq 0\right) \right) dq \\
&&\,\sim g_{\alpha ,\beta }(0)b_{n}a_{r}\int_{-z}^{(y-z)\wedge
0}V^{-}(-qa_{r})V^{+}\left( \left( t-z-q\right) a_{r}\right) dq.
\end{eqnarray*}%
Thus,
\begin{eqnarray*}
&&\mathbf{P}_{w}\left( S_{\tau _{r,n}}\leq ya_{r},\mathcal{B}(ta_{r},n),\tau
_{r,n}\in \lbrack n-\sqrt{nr}+r+1,n]\right) \\
&&\quad \sim g_{\alpha ,\beta }(0)b_{n}a_{r}\mathbf{P}\left( \tau
_{1}^{-}>r\right) V^{-}(w) \\
&&\qquad \times \int_{\varepsilon }^{N}g^{+}\left( z\right)
dz\int_{-z}^{(y-z)\wedge 0}V^{-}(-qa_{r})V^{+}\left( \left( t-z-q\right)
a_{r}\right) dq \\
&&\qquad =g_{\alpha ,\beta }(0)b_{n}\mathbf{P}\left( \tau _{1}^{-}>r\right)
a_{r}V^{-}(a_{r})V^{+}(a_{r})V^{-}(w) \\
&&\qquad \times \int_{\varepsilon }^{N}g^{+}\left( z\right)
dz\int_{-z}^{(y-z)\wedge 0}\frac{V^{-}(-qa_{r})V^{+}\left( \left(
t-z-q\right) a_{r}\right) }{V^{-}(a_{r})V^{+}(a_{r})}dq.
\end{eqnarray*}%
It follows from (\ref{Regular1}) and (\ref{Regular2}) that, as $r\rightarrow
\infty $
\begin{equation*}
\int_{-z}^{(y-z)\wedge 0}\frac{V^{-}(-qa_{r})V^{+}\left( \left( t-z-q\right)
a_{r}\right) }{V^{-}(a_{r})V^{+}(a_{r})}dq\sim \int_{\left( z-y\right) \vee
0}^{z}q^{\alpha \rho }\left( t-z+q\right) ^{\alpha (1-\rho )}dq
\end{equation*}%
uniformly in $z$ from any bounded interval $[0,N]$. Hence, using (\ref%
{Regular2}), (\ref{Product11}) and the equivalence%
\begin{equation*}
a_{r}V^{+}(a_{r})\sim (\alpha \rho +1)t^{-\alpha \rho
-1}\int_{0}^{ta_{r}}V^{+}(u)du, \, r\rightarrow \infty,
\end{equation*}%
which follows from (\ref{AsympV}) and (\ref{Regular1}), we conclude that, as
$n\gg r\rightarrow \infty $ and $\varepsilon \rightarrow 0$ and $%
N\rightarrow \infty $
\begin{eqnarray}
&&\mathbf{P}_{w}\left( S_{\tau _{r,n}}\leq ya_{r},\mathcal{B}(ta_{r},n),\tau
_{r,n}\in \lbrack n-\sqrt{nr},n]\right)  \notag \\
&\sim &V^{-}(w)g_{\alpha ,-\beta }(0)b_{n}\mathbf{P}\left( \tau
_{1}^{-}>r\right) V^{-}(a_{r})a_{r}V^{+}(a_{r})  \notag \\
&&\times \int_{\varepsilon }^{N}g^{+}\left( z\right) dz\int_{\left(
z-y\right) \vee 0}^{z}q^{\alpha \rho }\left( t-z+q\right) ^{\alpha (1-\rho
)}dq  \notag \\
&\sim &\frac{C^{\ast \ast }(\alpha \rho +1)}{t^{\alpha \rho +1}}%
V^{-}(w)g_{\alpha ,-\beta }(0)b_{n}\int_{0}^{ta_{r}}V^{+}(u)du  \notag \\
&&\times \int_{\varepsilon }^{N}g^{+}\left( z\right) dz\int_{\left(
z-y\right) \vee 0}^{z}q^{\alpha \rho }\left( t-z+q\right) ^{\alpha (1-\rho
)}dq  \notag \\
&\sim &\frac{C^{\ast \ast }(\alpha \rho +1)}{t^{\alpha \rho +1}}%
\int_{0}^{\infty }g^{+}\left( z\right) dz\int_{\left( z-y\right) \vee
0}^{z}q^{\alpha \rho }\left( t-z+q\right) ^{\alpha (1-\rho )}dq  \notag \\
&&\times \mathbf{P}_{w}\left( S_{\tau _{r,n}}\leq ya_{r},\mathcal{B}%
(ta_{r},n)\right) .  \label{Contribution1}
\end{eqnarray}

We now investigate the asymptotic behavior of $Q_{z}\left( 0,\sqrt{nr}%
\right) $. It is not difficult to check using (\ref{BasicAsymptotic}), (\ref%
{Defb}), and (\ref{Regular1}) that if $j\sim n$ then
\begin{equation*}
\mathbf{P}\left( S_{j}\leq \left( t-z\right) a_{r}-s,L_{j}\geq 0\right) \sim
\mathbf{P}\left( S_{n}\leq \left( t-z\right) a_{r}-s,L_{n}\geq 0\right)
\end{equation*}%
as $n\gg r\rightarrow \infty $ uniformly in $s\in \lbrack
-za_{r},((y-z)\wedge 0)a_{r})$. Therefore, as $n\gg r\rightarrow \infty $
\begin{eqnarray*}
&&Q_{z}\left( 0,\sqrt{nr}\right) \\
&=&\sum_{j=n-\sqrt{nr}+1}^{n-r}\int_{-za_{r}}^{((y-z)\wedge 0)a_{r}}\mathbf{P%
}\left( S_{n-r-j}\in ds,\tau _{n-r-j}=n-r-j\right) \mathbf{P}\left(
S_{j}\leq \left( t-z\right) a_{r}-s,L_{j}\geq 0\right) \\
&\sim &\int_{-z}^{(y-z)\wedge 0}\mathbf{P}\left( S_{n}\leq \left(
t-z-q\right) a_{r},L_{n}\geq 0\right) \sum_{l=0}^{\sqrt{nr}-r-1}\mathbf{P}%
\left( S_{l}\in a_{r}dq,\tau _{l}=l\right) \\
&=&\int_{-z}^{(y-z)\wedge 0}\mathbf{P}\left( S_{n}\leq \left( t-z-q\right)
a_{r},L_{n}\geq 0\right) V^{-}(-a_{r}dq)-\bar{Q}_{z}\left( 0,\sqrt{nr}%
\right) ,
\end{eqnarray*}%
where%
\begin{eqnarray*}
&&\bar{Q}_{z}\left( 0,\sqrt{nr}\right) \\
&&\quad =\int_{(z-y)\vee 0}^{z}\mathbf{P}\left( S_{n}\leq \left(
t-z+q\right) a_{r},L_{n}\geq 0\right) \sum_{l=\sqrt{nr}-r}^{\infty }\mathbf{P%
}\left( S_{l}\in -a_{r}dq,\tau _{l}=l\right) \\
&&\quad \leq \mathbf{P}\left( S_{n}\leq ta_{r},L_{n}\geq 0\right) \sum_{l=%
\sqrt{nr}-r}^{\infty }\mathbf{P}\left( S_{l}\in \lbrack -za_{r},((y-z)\wedge
0)a_{r}),\tau _{l}=l\right) .
\end{eqnarray*}

We first analyse the behavior of the integral%
\begin{eqnarray*}
&&\int_{-z}^{(y-z)\wedge 0}\mathbf{P}\left( S_{n}\leq \left( t-q\right)
a_{r},L_{n}\geq 0\right) V^{-}(-a_{r}dq) \\
&=&\mathbf{P}\left( S_{n}\leq ta_{r},L_{n}\geq 0\right) \\
&&\times V^{-}(a_{r})\int_{(z-y)\vee 0}^{z}\frac{\mathbf{P}\left( S_{n}\leq
\left( t-z+q\right) a_{r},L_{n}\geq 0\right) }{\mathbf{P}\left( S_{n}\leq
ta_{r},L_{n}\geq 0\right) }\frac{V^{-}(a_{r}dq)}{V^{-}(a_{r})}
\end{eqnarray*}%
as $n>>r\rightarrow \infty .$ In view of (\ref{BasicAsymptotic}), (\ref%
{AsympV}) and properties of regularly varying functions
\begin{eqnarray*}
&&\int_{(z-y)\vee 0}^{z}\frac{\mathbf{P}\left( S_{n}\leq \left( t-z+q\right)
a_{r},L_{n}\geq 0\right) }{\mathbf{P}\left( S_{n}\leq ta_{r},L_{n}\geq
0\right) }\frac{V^{-}(a_{r}dq)}{V^{-}(a_{r})} \\
&&\qquad \qquad \sim t^{-\alpha \rho -1}\int_{(z-y)\vee 0}^{z}\left(
t-z+q\right) ^{\alpha \rho +1}\frac{V^{-}(a_{r}dq)}{V^{-}(a_{r})}
\end{eqnarray*}%
uniformly in $z\in \lbrack \varepsilon ,N]$ as $n\gg r\rightarrow \infty $.
Since
\begin{equation*}
\frac{V^{-}(a_{r}q)}{V^{-}(a_{r})}\rightarrow q^{\alpha (1-\rho )},0\leq
q<\infty ,
\end{equation*}%
as $r\rightarrow \infty $, the measure%
\begin{equation*}
\frac{V^{-}(a_{r}q)}{V^{-}(a_{r})},0\leq q<\infty ,
\end{equation*}%
converges to the measure with density $\alpha (1-\rho )q^{\alpha (1-\rho
)-1} $. Thus, as $r\rightarrow \infty $

\begin{eqnarray*}
&&\int_{(z-y)\vee 0}^{z}\frac{\mathbf{P}\left( S_{n}\leq \left( t-z+q\right)
a_{r},L_{n}\geq 0\right) }{\mathbf{P}\left( S_{n}\leq ta_{r},L_{n}\geq
0\right) }\frac{V^{-}(a_{r}dq)}{V^{-}(a_{r})} \\
&\rightarrow &\alpha (1-\rho )t^{-\alpha \rho -1}\int_{(z-y)\vee
0}^{z}\left( t-z+q\right) ^{\alpha \rho +1}q^{\alpha (1-\rho )-1}dq.
\end{eqnarray*}

Our next goal is to evaluate $\bar{Q}_{z}\left( 0,\sqrt{nr}\right) $ from
above. By\ the duality principle for random walks, point (4) in Lemma~\ref%
{L_Renewal_negl} and the estimate (\ref{B_sum}) we have
\begin{eqnarray*}
&&\sum_{l=\sqrt{nr}-r}^{\infty }\mathbf{P}\left( S_{l}\in \lbrack
-za_{r},((y-z)\wedge 0)a_{r}),\tau _{l}=l\right) \\
&&\qquad \leq C_{1}\int_{(z-y)a_{r}\vee 0}^{za_{r}}V^{-}(s)ds\sum_{l=\sqrt{nr%
}-r}^{\infty }b_{j} \\
&&\qquad \leq C_{2}\frac{a_{r}V^{-}(a_{r})}{a_{\sqrt{nr}-r}}\int_{(z-y)\vee
0}^{z}q^{\alpha (1-\rho )}dq\leq C_{3}N^{\alpha (1-\rho )+1}\frac{%
a_{r}V^{-}(a_{r})}{a_{\sqrt{nr}}}.
\end{eqnarray*}

Thus,%
\begin{equation*}
\bar{Q}_{z}\left( 0,\sqrt{nr}\right) \leq C_{4}\mathbf{P}\left( S_{n}\leq
ta_{r},L_{n}\geq 0\right) N^{\alpha (1-\rho )+1}\frac{a_{r}V^{-}(a_{r})}{a_{%
\sqrt{nr}}}
\end{equation*}%
and, therefore, for fixed $0<\varepsilon <N<\infty $%
\begin{eqnarray*}
&&\mathbf{P}\left( \tau _{1}^{-}>r\right) V^{-}(w)\int_{\varepsilon
}^{N}\left( g^{+}\left( z\right) +o(1)\right) \bar{Q}\left( 0,\sqrt{nr}%
\right) dz \\
&\leq &C_{5}\frac{a_{r}V^{-}(a_{r})\mathbf{P}\left( \tau ^{-}>r\right) }{a_{%
\sqrt{nr}}}\leq C_{6}\frac{a_{r}}{a_{\sqrt{nr}}}\rightarrow 0
\end{eqnarray*}%
as $r\rightarrow \infty $. As a result, the following sequence of estimates
takes place:
\begin{eqnarray}
&&\mathbf{P}_{w}\left( S_{\tau _{r,n}}\leq ya_{r},\mathcal{B}(ta_{r},n),\tau
_{r,n}\in \lbrack r,\sqrt{nr}+r]\right)  \notag \\
&\sim &\mathbf{P}\left( \tau _{1}^{-}>r\right) V^{-}(w)\int_{\varepsilon
}^{N}g^{+}\left( z\right) \mathbf{P}_{za_{r}}\left( 0\leq S_{\tau
_{n-r}}\leq \left( y\wedge z\right) a_{r},S_{n-r}\leq ta_{r}\right) dz
\notag \\
&\sim &\mathbf{P}\left( \tau _{1}^{-}>r\right) V^{-}(w)  \notag \\
&&\quad \times \int_{\varepsilon }^{N}g^{+}\left( z\right)
\int_{-z}^{(y-z)\wedge 0}\mathbf{P}\left( S_{n}\leq \left( t-z-q\right)
a_{r},L_{n}\geq 0\right) V^{-}(-a_{r}dq)dz  \notag \\
&\sim &\mathbf{P}\left( \tau _{1}^{-}>r\right) V^{-}(a_{r})\alpha (1-\rho
)t^{-\alpha \rho -1}V^{-}(w)\mathbf{P}\left( S_{n}\leq ta_{r},L_{n}\geq
0\right) \times  \notag \\
&&\times \int_{\varepsilon }^{N}g^{+}\left( z\right) \int_{(z-y)\vee
0}^{z}\left( t-z+q\right) ^{\alpha \rho +1}q^{\alpha (1-\rho )-1}dqdz  \notag
\\
&\sim &C^{\ast \ast }\alpha (1-\rho )t^{-\alpha \rho -1}\int_{0}^{\infty
}g^{+}\left( z\right) \int_{(z-y)\vee 0}^{z}\left( t-z+q\right) ^{\alpha
\rho +1}q^{\alpha (1-\rho )-1}dqdz  \notag \\
&&\times \mathbf{P}_{w}\left( S_{n}\leq ta_{r},L_{n}\geq 0\right) ,
\label{Contribution2}
\end{eqnarray}%
where the symbol $\sim $ should be understood as the equivalence as $%
\varepsilon \downarrow 0,N\uparrow \infty ,$ and $n>>r\rightarrow \infty .$

Finally, we evaluate the intermediate term $Q(\sqrt{nr}+1,n-\sqrt{nr}).$ In
view of (\ref{Regular1}) and (\ref{Regular2})
\begin{eqnarray*}
\int_{0}^{\left( t+z\right) a_{r}}V^{+}(u)du &=&\frac{\int_{0}^{\left(
t+z\right) a_{r}}V^{+}(u)du}{\int_{0}^{ta_{r}}V^{+}(u)du}\times
\int_{0}^{ta_{r}}V^{+}(u)du \\
&\leq &C\left( t+z\right) ^{\alpha \rho +1}\int_{0}^{ta_{r}}V^{+}(u)du,
\end{eqnarray*}%
and, by (\ref{IntegralV_minus})
\begin{eqnarray*}
\int_{\left( (z-y)\wedge 0\right) a_{r}}^{za_{r}}V^{-}(q)dq &\leq
&a_{r}\int_{0}^{z}V^{-}(la_{r})dl \\
&\leq &C_{1}a_{r}V^{-}(a_{r})\int_{0}^{z}l^{\alpha (1-\rho
)}dl=C_{2}a_{r}V^{-}(a_{r})z^{\alpha (1-\rho )+1}.
\end{eqnarray*}%
Using the estimates above and (\ref{L_estimB}) we get%
\begin{eqnarray*}
&&Q_{z}(\sqrt{nr}+1,n-\sqrt{nr}) \\
&=&\sum_{j=\sqrt{nr}+1}^{n-\sqrt{nr}}\int_{-za_{r}}^{((y-z)\wedge 0)a_{r}}%
\mathbf{P}\left( S_{n-r-j}\in ds,\tau _{n-r-j}=n-r-j\right) \mathbf{P}\left(
S_{j}\leq \left( t-z\right) a_{r}-s,L_{j}\geq 0\right) \\
&\leq &\sum_{j=\sqrt{nr}+1}^{n-\sqrt{nr}}\mathbf{P}\left( S_{j}\leq
ta_{r},L_{j}\geq 0\right) \\
&&\times \mathbf{P}\left( S_{n-r-j}\in \lbrack -za_{r},((y-z)\wedge
0)a_{r}),\tau _{n-r-j}=n-r-j\right) \\
&\leq &C\sum_{j=\sqrt{nr}+1}^{n-\sqrt{nr}}b_{j}b_{n-r-j}%
\int_{0}^{ta_{r}}V^{+}(u)du\int_{\left( (z-y)\vee 0\right)
a_{r}}^{za_{r}}V^{-}(q)dq \\
&\leq &C_{5}\frac{a_{r}}{a_{\sqrt{nr}}}V^{-}(a_{r})z^{\alpha (1-\rho
)+1}b_{n}\int_{0}^{ta_{r}}V^{+}(u)du \\
&\leq &C_{6}\frac{a_{r}}{a_{\sqrt{nr}}}V^{-}(a_{r})N^{\alpha (1-\rho )+1}%
\mathbf{P}\left( \mathcal{B}(ta_{r},n)\right) .
\end{eqnarray*}%
Hence it follows that, for sufficiently large $r$ and $n\geq r$%
\begin{eqnarray*}
&&\mathbf{P}_{w}\left( S_{\tau _{r,n}}\leq ya_{r},\mathcal{B}(ta_{r},n),\tau
_{r,n}\in \lbrack \sqrt{nr}+1,n-\sqrt{nr}]\right) \\
&\leq &2\mathbf{P}\left( \tau _{1}^{-}>r\right) V^{-}(w)\int_{\varepsilon
}^{N}g^{+}\left( z\right) Q_{z}(\sqrt{nr}+1,n-\sqrt{nr})dz \\
&&+\mathbf{P}_{w}\left( S_{\tau _{r,n}}\leq ya_{r},S_{r}\notin \lbrack
\varepsilon a_{r},Na_{r}],\mathcal{B}(ta_{r},n),\tau _{r,n}\in \lbrack \sqrt{%
nr}+1,n-\sqrt{nr}]\right) \\
&\leq &C_{3}\mathbf{P}\left( \mathcal{B}(ta_{r},n)\right) \frac{a_{r}}{a_{%
\sqrt{nr}}}V^{-}(a_{r})\mathbf{P}\left( \tau _{1}^{-}>r\right)
V^{-}(w)\int_{\varepsilon }^{N}g^{+}\left( z\right) dzN^{\alpha (1-\rho )+1}
\\
&&+\mathbf{P}_{w}\left( S_{r}\notin \lbrack \varepsilon a_{r},Na_{r}],%
\mathcal{B}(ta_{r},n)\right) \\
&\leq &2C^{\ast \ast }C_{3}V^{-}(w)\mathbf{P}\left( \mathcal{B}%
(ta_{r},n)\right) \frac{a_{r}}{a_{\sqrt{nr}}}N^{\alpha (1-\rho )+1} \\
&&+\mathbf{P}_{w}\left( S_{r}\notin \lbrack \varepsilon a_{r},Na_{r}],%
\mathcal{B}(ta_{r},n)\right) ,
\end{eqnarray*}%
where
\begin{equation*}
\frac{a_{r}}{a_{\sqrt{nr}}}N^{\alpha (1-\rho )+1}\rightarrow 0
\end{equation*}%
as $n\gg r\rightarrow \infty $. Therefore,
\begin{equation*}
\mathbf{P}_{w}\left( S_{\tau _{r,n}}\leq ya_{r},\mathcal{B}(ta_{r},n),\tau
_{r,n}\in \lbrack \sqrt{nr}+1,n-\sqrt{nr}]\right) =o\left( \mathbf{P}\left(
\mathcal{B}(ta_{r},n)\right) \right)
\end{equation*}%
as $n\gg r\rightarrow \infty $. Combining this estimate with (\ref%
{Contribution1})--(\ref{Contribution2}) we see that%
\begin{equation*}
\mathbf{P}_{w}\left( S_{\tau _{r,n}}\leq ya_{r},\mathcal{B}(ta_{r},n)\right)
\sim W(t,y)\mathbf{P}_{w}\left( \mathcal{B}(ta_{r},n)\right) ,
\end{equation*}%
as $n\gg r\rightarrow \infty $, where%
\begin{eqnarray*}
W(t,y) &=&\frac{C^{\ast \ast }(\alpha \rho +1)}{t^{\alpha \rho +1}}%
\int_{0}^{\infty }g^{+}\left( z\right) dz\int_{(z-y)\vee 0}^{z}q^{\alpha
\rho }(t-z+q)^{\alpha (1-\rho )}dq \\
&&+\frac{C^{\ast \ast }\alpha (1-\rho )}{t^{\alpha \rho +1}}\int_{0}^{\infty
}g^{+}\left( z\right) dz\int_{(z-y)\vee 0}^{z}\left( t-z+q\right) ^{\alpha
\rho +1}q^{\alpha (1-\rho )-1}dq.
\end{eqnarray*}

Using again (\ref{Proport0}),(\ref{Inclusion}) and (\ref{BasicAsymptotic}) \
we conclude that
\begin{eqnarray*}
\mathbf{P}_{w}\left( S_{\tau _{r,n}}\leq ya_{r},\mathcal{B}(ta_{k},n)\right)
&\sim &\mathbf{P}_{w}\left( S_{\tau _{r,n}}\leq ya_{r},\mathcal{B}(t\theta
^{1/\alpha }a_{r},n)\right) \\
&\sim &W(t\theta ^{1/\alpha },y)\mathbf{P}_{w}\left( \mathcal{B}(t\theta
^{1/\alpha }a_{r},n)\right) \\
&\sim &W(t\theta ^{1/\alpha },y)\mathbf{P}_{w}\left( \mathcal{B}%
(ta_{k},n)\right)
\end{eqnarray*}%
as \ $n\gg k=[\theta r]\rightarrow \infty $.

Theorem \ref{T_initialDimilar} is proved.

\section{The case $\min \left( r,n-r\right) \gg k$\label{Sec4}}

First we assume that $n\gg r\gg k\rightarrow \infty $ and prove the
following statement.

\begin{theorem}
\label{T_veryFar} If Conditions A1 and A2 are valid, then, for any $w\geq
0,t\in \left( 0,\infty \right) $ and $y\in \lbrack 0,t]$
\begin{equation*}
\lim_{n\gg r\gg k\rightarrow \infty }\mathbf{P}_{w}\left( S_{\tau
_{r,n}}\leq ya_{k}|\mathcal{B}(ta_{k},n)\right) =1-\left( 1-\frac{y}{t}%
\right) ^{\alpha \rho +1}.
\end{equation*}
\end{theorem}

\textbf{Proof}. The same as in Theorem \ref{T_beginningNew} we fix $%
0<\varepsilon <N<\infty $ and use the representation
\begin{eqnarray*}
&&\mathbf{P}_{w}\left( S_{\tau _{r,n}}\leq ya_{k},\mathcal{B}%
(ta_{k},n)\right) \\
&\sim &\mathbf{P}\left( \tau _{1}^{-}>r\right) V^{-}(w)\int_{\varepsilon
}^{N}g^{+}\left( z\right) \mathbf{P}_{za_{r}}\left( 0\leq S_{\tau
_{n-r}}\leq ya_{k},S_{n-r}\leq ta_{k}\right) dz \\
&&+\mathbf{P}_{w}\left( S_{\tau _{r,n}}\leq ya_{k},S_{r}\notin \lbrack
\varepsilon a_{r},Na_{r}],\mathcal{B}(ta_{k},n)\right) ,
\end{eqnarray*}%
where
\begin{equation}
\lim_{\varepsilon \downarrow 0,N\uparrow \infty }\lim_{n\gg r\gg
k\rightarrow \infty }\frac{\mathbf{P}_{w}\left( S_{\tau _{r,n}}\leq
ya_{k},S_{r}\notin \lbrack \varepsilon a_{r},Na_{r}],\mathcal{B}%
(ta_{k},n)\right) }{\mathbf{P}_{w}\left( \mathcal{B}(ta_{k},n)\right) }=0
\label{VeryFar1}
\end{equation}%
according to Corollary \ref{C_conditNegligible}. We select a sequence $%
m=m(n) $ such that $r\gg m\gg k$ and write
\begin{eqnarray}
&&\mathbf{P}_{a_{r}z}\left( 0\leq S_{\tau _{n-r}}\leq ya_{k},S_{n-r}\leq
ta_{r}\right) =Q_{z}(0,\sqrt{mk})  \notag \\
&&+Q_{z}(\sqrt{mk}+1,n-r-\sqrt{mk})+Q_{z}(n-r-\sqrt{mk}+1,n-r),
\label{Decompos1}
\end{eqnarray}%
where now%
\begin{equation*}
Q_{z}(A,B)=\mathbf{P}_{za_{r}}\left( 0\leq S_{\tau _{n-r}}\leq
ya_{k},S_{n-r}\leq ta_{r};\tau _{n-r}\in \lbrack A,B]\right) .
\end{equation*}%
First we estimate the term
\begin{eqnarray*}
&&Q_{z}(n-r-\sqrt{mk}+1,n-r) \\
&&\qquad =\sum_{j=0}^{\sqrt{mk}-1}\int_{-za_{r}}^{ya_{k}-za_{r}}\mathbf{P}%
\left( S_{n-r-j}\in ds,\tau _{n-r-j}=n-r-j\right) \\
&&\qquad \quad \times \mathbf{P}\left( S_{j}\leq \left( t-z\right)
a_{k}-s,L_{j}\geq 0\right) .
\end{eqnarray*}%
By the duality principle for random walks and formula (2.17) in \cite%
{VDD2023}, applied to $-S,$ we have for $n\gg r\gg k$ and $j\in \lbrack 0,%
\sqrt{mk}]:$
\begin{eqnarray*}
&&\mathbf{P}\left( S_{n-r-j}\in ds,\tau _{n-r-j}=n-r-j\right) =\mathbf{P}%
\left( S_{n-r-j}\in ds,M_{n-r-j}<0\right) \\
&&\qquad \qquad \qquad \sim g_{\alpha ,\beta }(0)b_{n-j-r}V^{-}(-s)ds\sim
g_{\alpha ,\beta }(0)b_{n}V^{-}(za_{r})ds.
\end{eqnarray*}%
Using this estimate and relation (\ref{Tail11}), we conclude that
\begin{eqnarray*}
&&Q_{z}(n-r-\sqrt{mk}+1,n-r) \\
&&\quad \sim g_{\alpha ,\beta
}(0)b_{n}V^{-}(za_{r})\int_{-za_{r}}^{ya_{k}-za_{r}}\sum_{j=0}^{\sqrt{mk}}%
\mathbf{P}\left( S_{j}\leq \left( t-z\right) a_{k}-s,L_{j}\geq 0\right) ds \\
&&\qquad =g_{\alpha ,\beta }(0)b_{n}V^{-}(za_{r})\text{ }\int_{0}^{ya_{k}}%
\sum_{j=0}^{\sqrt{mk}}\mathbf{P}\left( S_{j}\leq ta_{k}-s,L_{j}\geq 0\right)
ds \\
&&\qquad \sim g_{\alpha ,\beta }(0)b_{n}V^{-}(za_{r})\text{ }%
\int_{0}^{ya_{k}}V^{+}\left( ta_{k}-s\right) ds \\
&&\qquad =g_{\alpha ,\beta
}(0)b_{n}V^{-}(za_{r})\int_{(t-y)a_{k}}^{ta_{k}}V^{+}(s)ds
\end{eqnarray*}%
as $n\gg r\gg k\rightarrow \infty .$ Therefore,
\begin{eqnarray*}
&&\int_{\varepsilon }^{N}\mathbf{P}_{w}\left( S_{r}\in a_{r}dz;L_{r}\geq
0\right) Q_{z}(n-r-\sqrt{mk}+1,n-r) \\
&\sim &g_{\alpha ,\beta }(0)V^{-}(w)b_{n}\int_{(t-y)a_{k}}^{ta_{k}}V^{+}(s)ds%
\mathbf{P}\left( \tau _{1}^{-}>r\right) \int_{\varepsilon }^{N}g^{+}\left(
z\right) V^{-}(za_{r})dz.
\end{eqnarray*}%
By (\ref{Regular2}) and (\ref{IntegralV_plus})
\begin{eqnarray*}
\mathbf{P}\left( \tau _{1}^{-}>r\right) \int_{\varepsilon }^{N}g^{+}\left(
z\right) V^{-}(a_{r}z)dz &\sim &\mathbf{P}\left( \tau _{1}^{-}>r\right)
V^{-}(a_{r})\int_{\varepsilon }^{N}g^{+}\left( z\right) z^{\alpha (1-\rho
)}dz \\
&\sim &C^{\ast \ast }\int_{\varepsilon }^{N}g^{+}\left( z\right) z^{\alpha
(1-\rho )}dz
\end{eqnarray*}%
as $r\rightarrow \infty $. Thus,
\begin{eqnarray}
&&\mathbf{P}_{w}\left( S_{\tau _{r,n}}\leq ya_{k},\mathcal{B}(ta_{k},n),\tau
_{r,n}\in \lbrack n-\sqrt{mk}+1,n]\right)  \notag \\
&=&(1+o(1))\int_{\varepsilon }^{N}\mathbf{P}_{w}\left( S_{r}\in
a_{r}dz;L_{r}\geq 0\right) Q_{z}(n-r-\sqrt{mk}+1,n-r)  \notag \\
&\sim &C^{\ast \ast }g_{\alpha ,\beta }(0)V^{-}(w)b_{n}\text{ }%
\int_{(t-y)a_{k}}^{ta_{k}}V^{+}(s)ds\int_{\varepsilon }^{N}g^{+}\left(
z\right) z^{\alpha (1-\rho )}dz.  \label{VeryFar2}
\end{eqnarray}%
We now evaluate the term
\begin{equation*}
Q_{z}(0,\sqrt{mk})=\sum_{j=0}^{\sqrt{mk}}\int_{0}^{ya_{k}}\mathbf{P}%
_{za_{r}}\left( S_{j}\in ds,\tau _{j}=j\right) \mathbf{P}\left(
S_{n-r-j}\leq ta_{k}-s,L_{n-r-j}\geq 0\right) .
\end{equation*}%
In view of (\ref{BasicAsymptotic}) and properties of regularly varying
functions the estimate
\begin{eqnarray*}
\mathbf{P}\left( S_{n-r-j}\leq ta_{k}-s,L_{n-r-j}\geq 0\right) &\leq &%
\mathbf{P}\left( S_{n-r-j}\leq ta_{k},L_{n-r-j}\geq 0\right) \\
&\leq &C\mathbf{P}\left( S_{n}\leq ta_{k},L_{n}\geq 0\right)
\end{eqnarray*}%
is valid for $j\leq \sqrt{mk}$ as $n\gg r\gg k\rightarrow \infty .$ Further,
we have
\begin{eqnarray*}
\sum_{j=0}^{\sqrt{mk}}\mathbf{P}_{za_{r}}\left( S_{j}\in \lbrack
0,ya_{k}),\tau _{j}=j\right) &\leq &\sum_{j=0}^{\infty }\mathbf{P}%
_{za_{r}}\left( S_{j}\in \lbrack 0,ya_{k}),\tau _{j}=j\right) \\
&=&\hat{V}^{-}(za_{r})-\hat{V}^{-}(za_{r}-ya_{k}).
\end{eqnarray*}%
Thus,
\begin{equation*}
Q_{z}(0,\sqrt{mk})\leq C\mathbf{P}\left( \mathcal{B}(ta_{k},n)\right) (\hat{V%
}^{-}(za_{r})-\hat{V}^{-}(za_{r}-ya_{k})).
\end{equation*}%
Since $\mathbf{P}\left( \tau _{1}^{-}>r\right) V^{-}(a_{r})\sim C^{\ast \ast
}$ as $r\rightarrow \infty $ in view of (\ref{Product11}) and
\begin{equation*}
\frac{\hat{V}^{-}(za_{r})-\hat{V}^{-}(za_{r}-ya_{k})}{V^{-}(a_{r})}=(1-\zeta
)\frac{\hat{V}^{-}(za_{r})-\hat{V}^{-}(za_{r}-ya_{k})}{\hat{V}^{-}(a_{r})}%
\rightarrow 0
\end{equation*}%
uniformly in $z\in \lbrack \varepsilon ,N]$ as $r\gg k\rightarrow \infty $,
we conclude that
\begin{eqnarray}
&&\mathbf{P}\left( \tau _{1}^{-}>r\right) V^{-}(w)\int_{\varepsilon
}^{N}g^{+}\left( z\right) Q_{z}(0,\sqrt{mk})dz  \notag \\
&&\qquad \leq CV^{-}(w)\mathbf{P}\left( \mathcal{B}(ta_{k},n)\right) \mathbf{%
P}\left( \tau _{1}^{-}>r\right) V^{-}(a_{r})  \notag \\
&&\qquad \times \int_{\varepsilon }^{N}g^{+}\left( z\right) \frac{\hat{V}%
^{-}(za_{r})-\hat{V}^{-}(za_{r}-ya_{k})}{V^{-}(a_{r})}dz  \notag \\
&&\qquad \qquad =o\left( \mathbf{P}\left( \mathcal{B}(ta_{k},n)\right)
\right)  \label{VeryFar3}
\end{eqnarray}%
as $n\gg r\gg k\rightarrow \infty $. Now we estimate the term
\begin{eqnarray*}
&&Q_{z}(\sqrt{mk}+1,n-r-\sqrt{mk}) \\
&&\quad =\sum_{j=\sqrt{mk}+1}^{n-r-\sqrt{mk}}\int_{0}^{ya_{k}}\mathbf{P}%
_{za_{r}}\left( S_{n-r-j}\in ds,\tau _{n-r-j}=n-r-j\right) \\
&&\qquad \qquad \times \mathbf{P}\left( S_{j}\leq ta_{k}-s,L_{j}\geq
0\right) .
\end{eqnarray*}%
Using points (3) and (4) of Lemma \ref{L_Renewal_negl} and Lemma \ref%
{L_estimB} we see that
\begin{eqnarray*}
&&Q_{z}(\sqrt{mk}+1,n-r-\sqrt{mk}) \\
&\leq &C\int_{0}^{ta_{k}}V^{+}(u)du \\
&&\qquad \times \sum_{j=\sqrt{mk}+1}^{n-r-\sqrt{mk}}b_{j}\mathbf{P}\left(
-za_{r}\leq S_{n-r-j}<-za_{r}+ya_{k},\tau _{n-r-j}=n-r-j\right) \\
&\leq
&C\int_{0}^{ta_{k}}V^{+}(u)du\int_{-za_{r}}^{-za_{r}+ya_{k}}V^{-}(q)dq%
\sum_{j=\sqrt{mk}+1}^{n-r-\sqrt{mk}}b_{j}b_{n-r-j} \\
&\leq &C_{1}y\frac{a_{k}}{a_{\sqrt{mk}}}b_{n}V^{-}(za_{r})%
\int_{0}^{ta_{k}}V^{+}(u)du.
\end{eqnarray*}%
Therefore,
\begin{eqnarray}
&&\mathbf{P}\left( \tau _{1}^{-}>r\right) V^{-}(w)\int_{\varepsilon
}^{N}g^{+}\left( z\right) Q_{z}(\sqrt{mk}+1,n-r-\sqrt{mk})dz  \notag \\
&&\,\leq C_{1}y\frac{a_{k}}{a_{\sqrt{mk}}}b_{n}\int_{0}^{ta_{k}}V^{+}(u)du\,%
\mathbf{P}\left( \tau _{1}^{-}>r\right) \int_{\varepsilon }^{N}g^{+}\left(
z\right) V^{-}(za_{r})dz  \notag \\
&&\,\leq C_{2}y\frac{a_{k}}{a_{\sqrt{mk}}}b_{n}\int_{0}^{ta_{k}}V^{+}(u)du\,%
\mathbf{P}\left( \tau _{1}^{-}>r\right) V^{-}(a_{r})\int_{\varepsilon
}^{N}g^{+}\left( z\right) z^{\alpha (1-\rho )}dz  \notag \\
&&\,\leq C_{3}y\frac{a_{k}}{a_{\sqrt{mk}}}b_{n}\int_{0}^{ta_{k}}V^{+}(u)du%
\int_{\varepsilon }^{N}g^{+}\left( z\right) z^{\alpha (1-\rho )}dz  \notag \\
&&\qquad =o\left( b_{n}\int_{0}^{ta_{k}}V^{+}(u)du\right) ,  \label{VeryFar4}
\end{eqnarray}%
since $a_{k}/a_{\sqrt{mk}}\rightarrow 0$ as $m\gg k\rightarrow \infty $.

Combining (\ref{VeryFar1})-(\ref{VeryFar4}), we obtain
\begin{eqnarray}
&&\mathbf{P}_{w}\left( S_{\tau _{r,n}}\leq ya_{k},\mathcal{B}%
(ta_{k},n)\right)  \notag \\
&\sim &\mathbf{P}_{w}\left( S_{\tau _{r,n}}\leq ya_{k},\mathcal{B}%
(ta_{k},n),\tau _{r,n}\in \lbrack n-\sqrt{mk}+1,n]\right)  \notag \\
&\sim &C^{\ast \ast }g_{\alpha ,\beta }(0)V^{-}(w)b_{n}\text{ }%
\int_{(t-y)a_{k}}^{ta_{k}}V^{+}(s)ds\int_{0}^{\infty }g^{+}\left( z\right)
z^{\alpha (1-\rho )}dz.  \label{MinimLocalization}
\end{eqnarray}%
Setting here $y=t$, we conclude that
\begin{eqnarray*}
&&\mathbf{P}_{w}\left( \mathcal{B}(ta_{k},n)\right) \sim \mathbf{P}%
_{w}\left( \mathcal{B}(ta_{k},n),\tau _{r,n}\in \lbrack n-\sqrt{mk}%
+1,n]\right) \\
&&\qquad \sim C^{\ast \ast }g_{\alpha ,\beta }(0)V^{-}(w)b_{n}\text{ }%
\int_{0}^{ta_{k}}V^{+}(s)ds\int_{0}^{\infty }g^{+}\left( z\right) z^{\alpha
(1-\rho )}dz.
\end{eqnarray*}%
Hence, observing that
\begin{equation*}
\int_{(t-y)a_{k}}^{ta_{k}}V^{+}(s)ds\big/\int_{0}^{ta_{k}}V^{+}(s)ds%
\rightarrow 1-\left( 1-\frac{y}{t}\right) ^{\alpha \rho +1}
\end{equation*}%
as $k\rightarrow \infty ,$ we deduce
\begin{equation*}
\lim_{n\gg r\gg k\rightarrow \infty }\mathbf{P}_{w}\left( S_{\tau
_{r,n}}\leq ya_{k}|\mathcal{B}(ta_{k},n)\right) =1-\left( 1-\frac{y}{t}%
\right) ^{\alpha \rho +1}.
\end{equation*}

Theorem \ref{T_veryFar} is proved.

\begin{corollary}
\label{C_minimum_localization}If Conditions A1 and A2 are valid then, for
fixed $w\geq 0,t\in (0,\infty )$ and any sequence $r=r(n)$ such that $\min
\left( r,n-r\right) \gg k\rightarrow \infty $%
\begin{equation}
\lim_{\min \left( r,n-r\right) \gg k\rightarrow \infty }\mathbf{P}_{w}\left(
S_{\tau _{r,n}}=S_{\tau _{n-r},n}|\mathcal{B}(ta_{k},n)\right) =1
\label{Global}
\end{equation}%
and, therefore,
\begin{equation}
\lim_{\min (r,n-r)\gg k\rightarrow \infty }\mathbf{P}_{w}\left( S_{\tau
_{r,n}}\leq ya_{k}|\mathcal{B}(ta_{k},n)\right) =1-\left( 1-\frac{y}{t}%
\right) ^{\alpha \rho +1}.  \label{Global2}
\end{equation}
\end{corollary}

\textbf{Proof}. Assume first that $r<n-r$. Take a sequence $m=m(n)$ such
that $r\gg m\gg k$. In view of (\ref{MinimLocalization}) \ \ \
\begin{eqnarray*}
\mathbf{P}_{w}\left( S_{\tau _{r,n}}\neq S_{\tau _{n-r},n},\mathcal{B}%
(ta_{k},n)\right) &=&\mathbf{P}_{w}\left( \tau _{r,n}\in \lbrack r,n-r),%
\mathcal{B}(ta_{k},n)\right) \\
&\leq &\mathbf{P}_{w}\left( \tau _{r,n}\in \lbrack r,n-\sqrt{mk}),\mathcal{B}%
(ta_{k},n)\right) \\
&=&o\left( \mathbf{P}_{w}\left( \mathcal{B}(ta_{k},n)\right) \right)
\end{eqnarray*}%
for $n\geq 2r\gg m\gg k$. The case $n-r<r$ may be considered in a similar
way. Relation (\ref{Global}) is proved.

Relation (\ref{Global2}) is an easy consequence of Theorem \ref{T_veryFar}.

\section{\protect\bigskip The case $n\gg k=\protect\theta \left( n-r\right) $%
\label{Sec5}}

We have considered in Section \ref{Sec4} the asymptotic behavior of the
probabilities
\begin{equation*}
\mathbf{P}_{w}\left( \frac{1}{a_{m}}S_{\tau _{r,n}}\leq y|\mathcal{B}%
(ta_{k},n)\right)
\end{equation*}%
for $m=n-r\gg k$. In this section we deal with the case when the difference $%
m=n-r$ is of order $k.$ As before, we agree the consider the values $\theta
\left( n-r\right) $ \ and $\theta m$ as \ $[\theta \left( n-r\right) ]$ \
and $[\theta m].$ We prove the following statement which is an extension of
Corollary 2 in \cite{VDD2023}.

\begin{theorem}
\label{T_interm} Let Conditions A1 and A2 be valid and $k\thicksim \theta
m\rightarrow \infty $ as $n\rightarrow \infty $. If $m=o\left( n\right) $ as
$n\rightarrow \infty ,$then, for any fixed $w\geq 0,t$ $\in (0,\infty )$ and
any $y\in \lbrack 0,t\theta ^{1/\alpha }]$
\begin{equation*}
\lim_{n\rightarrow \infty }\mathbf{P}_{w}\left( \frac{1}{a_{m}}S_{\tau
_{n-m,n}}\leq y|\mathcal{B}(Ta_{k},n)\right) =A(t\theta ^{1/\alpha },y),
\end{equation*}%
where%
\begin{equation*}
A(T,y):=\frac{\alpha \rho +1}{t^{1+\alpha \rho }}\int_{0}^{\infty }z^{\alpha
\rho }\mathbf{P}\left( -z\leq \min_{0\leq s\leq 1}Y_{s}\leq y-z,Y_{1}\leq
t-z\right) dz.
\end{equation*}
\end{theorem}

\textbf{Proof}. As in the proof of Theorem~\ref{T_initialDimilar} we first
find an asymptotic representation for the probability $\mathbf{P}_{w}\left(
S_{\tau _{n-m,n}}\leq ya_{m},\mathcal{B}(ta_{m},n)\right) ,$ and then use
analogs of the relations (\ref{Proport0}) and (\ref{Inclusion}). Set $%
x=ta_{m}$. Using (\ref{BasicAsymptotic}) and taking into account (\ref%
{AsympV}) we conclude that, for any $u\in (0,1)$
\begin{equation}
\mathbf{P}_{w}\mathbf{(}S_{n}\leq ux|\mathcal{B}(x,n))\rightarrow u^{\alpha
\rho +1}  \label{End_point_neglogible}
\end{equation}%
as $n\rightarrow \infty $ uniformly in \thinspace $x,w\geq 0$ such that $%
\max (x,w)\in (0,\delta _{n}a_{n}],$ where $\delta _{n}\rightarrow 0$ as $%
n\rightarrow \infty $. It is not difficult to show that for fixed positive
numbers $\varepsilon <N$
\begin{eqnarray*}
&&\mathbf{P}_{w}\left( \frac{1}{a_{m}}S_{\tau _{n-m,n}}\leq y,\mathcal{B}%
(x,n)\right) \\
&=&\mathbf{P}_{w}\left( \frac{1}{a_{m}}S_{\tau _{n-m,n}}\leq y,\varepsilon
a_{m}\leq S_{n-m}\leq Na_{m},\varepsilon x\leq S_{n}\leq x,L_{n}\geq 0\right)
\\
&&+r_{m,n}^{\ast }(\varepsilon ,N),
\end{eqnarray*}%
where
\begin{eqnarray*}
0\leq r_{m,n}^{\ast }(\varepsilon ,N) &\leq &\mathbf{P}_{w}\left( 0\leq
S_{n}\leq \varepsilon x,L_{n}\geq 0\right) \\
&&+\mathbf{P}_{w}\left( 0\leq S_{n-m}\leq \varepsilon a_{m},0\leq S_{n}\leq
x,L_{n}\geq 0\right) .
\end{eqnarray*}%
Observe that
\begin{equation}
\lim_{\varepsilon \downarrow 0,N\uparrow \infty }\limsup_{n\gg m\rightarrow
\infty }\frac{r_{m,n}^{\ast }(\varepsilon ,N)}{\mathbf{P}_{w}\mathbf{(}%
S_{n}\leq x,L_{n}\geq 0)}=0  \label{Interm1}
\end{equation}%
in view of (\ref{End_point_neglogible}) and Corollary 2 in \cite{VDD2023}.

For $x=ta_{m}$ introduce the event
\begin{equation*}
\mathcal{K}(m,n,\varepsilon ,x):=\{\varepsilon a_{m}\leq S_{n-m}\leq
Na_{m},\varepsilon x\leq S_{n}\leq x\}.
\end{equation*}%
Then
\begin{eqnarray*}
&&\mathbf{P}_{w}\left( \frac{1}{a_{m}}S_{\tau _{n-m,n}}\leq y,\varepsilon
a_{m}\leq S_{n-m}\leq Na_{m},\varepsilon x\leq S_{n}\leq x,L_{n}\geq 0\right)
\\
&=&\mathbf{P}_{w}\left( \frac{1}{a_{m}}(S_{\tau _{n-m,n}}-S_{n-m})\leq y-%
\frac{1}{a_{m}}S_{n-m},\mathcal{K}(m,n,\varepsilon ,x),L_{n}\geq 0\right) \\
&=&\int_{\varepsilon }^{N}\mathbf{P}_{w}\left( S_{n-m}\in
a_{m}dz,L_{n-m}\geq 0\right) \\
&&\times \mathbf{P}\left( -z\leq \frac{1}{a_{m}}S_{\tau _{m}}\leq
y-z,\varepsilon t-z\leq \frac{1}{a_{m}}S_{m}\leq t-z\right) .
\end{eqnarray*}%
Using the Donsker-Prokhorov invariance principle for random walks with
increments satisfying Condition A1, we obtain
\begin{eqnarray*}
&&\lim_{n\rightarrow \infty }\mathbf{P}\left( -z\leq \frac{1}{a_{m}}S_{\tau
_{m}}\leq y-z,\varepsilon t-z\leq \frac{1}{a_{m}}S_{m}\leq t-z\right) \\
&&\qquad =\mathbf{P}\left( -z\leq \min_{0\leq s\leq 1}Y_{s}\leq
y-z,\varepsilon t-z\leq Y_{1}\leq t-z\right) =:G(y,z,t,\varepsilon ).
\end{eqnarray*}%
By Theorem 5.1 in \cite{CC2013} and (\ref{Defb})
\begin{eqnarray}
\frac{\mathbf{P}_{w}\left( S_{n-m}\in a_{m}dz;L_{n-m}\geq 0\right) }{a_{m}dz}
&\sim &g_{\alpha ,\beta }(0)b_{n-m}V^{-}(w)V^{+}(a_{m}z)  \notag \\
&\sim &g_{\alpha ,\beta }(0)b_{n}V^{-}(w)V^{+}(a_{m}z)  \notag
\end{eqnarray}%
as $n\gg m\rightarrow \infty $ uniformly in $a_{m}z=o(a_{n})$. Using this
asymptotical representation and (\ref{IntegralV_plus}) we conclude that
\begin{eqnarray*}
&&\int_{\varepsilon }^{N}\mathbf{P}_{w}\left( S_{n-m}\in a_{m}dz,L_{n-m}\geq
0\right) \\
&&\qquad \times \mathbf{P}\left( -z\leq \frac{1}{a_{m}}S_{\tau _{m}}\leq
y-z,\varepsilon t-z\leq \frac{1}{a_{m}}S_{m}\leq t-z\right) \\
&\sim &g_{\alpha ,\beta }(0)b_{n-m}a_{m}V^{-}(w)\int_{\varepsilon
}^{N}V^{+}(a_{m}z)G(y,z,t,\varepsilon )dz \\
&\sim &g_{\alpha ,\beta }(0)b_{n}a_{m}V^{-}(w)V^{+}(a_{m})\int_{\varepsilon
}^{N}z^{\alpha \rho }G(y,z,t,\varepsilon )dz
\end{eqnarray*}%
as $m\rightarrow \infty $. By (\ref{Interm1}) we see that
\begin{eqnarray}
&&\mathbf{P}_{w}\left( \frac{1}{a_{m}}S_{\tau _{n-m,n}}\leq y,\mathcal{B}%
(ta_{m},n)\right)  \notag \\
&&\qquad \qquad \sim g_{\alpha ,\beta }(0)b_{n}a_{m}V^{-}(w)V^{+}(a_{m})
\notag \\
&&\qquad \qquad \times \int_{0}^{\infty }z^{\alpha \rho }\mathbf{P}\left(
-z\leq \min_{0\leq s\leq 1}Y_{s}\leq y-z,Y_{1}\leq t-z\right) dz
\label{Uncond1}
\end{eqnarray}%
as $n\gg m\rightarrow \infty $. Since
\begin{eqnarray*}
&&\mathbf{P}_{w}\left( \mathcal{B}(ta_{m},n)\right) \sim g_{\alpha ,\beta
}(0)b_{n}V^{-}(w)\int_{0}^{ta_{m}}V^{+}(z)dz \\
&&\quad \sim \frac{g_{\alpha ,\beta }(0)b_{n}}{\alpha \rho +1}%
V^{-}(w)ta_{m}V^{+}(ta_{m})\sim \frac{g_{\alpha ,\beta }(0)b_{n}}{\alpha
\rho +1}V^{-}(w)t^{1+\alpha \rho }a_{m}V^{+}(a_{m})
\end{eqnarray*}%
as $n\gg m\rightarrow \infty ,$ it follows that, for any $w\geq 0$ and $y\in
\lbrack 0,t]$%
\begin{eqnarray*}
&&\mathbf{P}_{w}\left( \frac{1}{a_{m}}S_{\tau _{n-m,n}}\leq y|\mathcal{B}%
(ta_{m},n)\right) \\
&&\qquad \sim \frac{\alpha \rho +1}{t^{1+\alpha \rho }}\int_{0}^{\infty
}z^{\alpha \rho }\mathbf{P}\left( -z\leq \min_{0\leq s\leq 1}Y_{s}\leq
y-z,Y_{1}\leq t-z\right) dz \\
&&\qquad \qquad =A(t,y)
\end{eqnarray*}%
as $n\gg m\rightarrow \infty $. To finish the proof of Theorem \ref{T_interm}
it remains to use (\ref{Proport0}) and (\ref{Inclusion}) for $r=m.$

\begin{remark}
Since the asymptotic representation (\ref{Uncond1}) is valid for $y=t$, we
may write
\begin{equation*}
A(t,y)=\frac{\int_{0}^{\infty }z^{\alpha \rho }\mathbf{P}\left( -z\leq
\min_{0\leq s\leq 1}Y_{s}\leq y-z,Y_{1}\leq t-z\right) dz}{\int_{0}^{\infty
}z^{\alpha \rho }\mathbf{P}\left( -z\leq \min_{0\leq s\leq 1}Y_{s},Y_{1}\leq
t-z\right) dz}.
\end{equation*}%
Taking into account that $A(t,t)=1$, we get an interesting relation
\begin{equation*}
\int_{0}^{\infty }z^{\alpha \rho }\mathbf{P}\left( -z\leq \min_{0\leq s\leq
1}Y_{s},Y_{1}\leq t-z\right) dz=\frac{t^{1+\alpha \rho }}{\alpha \rho +1}.
\end{equation*}
\end{remark}

\section{The case $n\gg k\gg n-r$\label{Sec6}}

In this section, to avoid cumbersome formulas we set $m=n-r$ and consider
the case $n\gg k\gg m\rightarrow \infty. $

\begin{theorem}
\label{T_Minim} Let Conditions A1, A2 be valid and $n\gg k\gg
m=n-r\rightarrow \infty .$ Then, for any $w\geq 0,t\in (0,\infty )$ and $%
y\leq 0$
\begin{equation*}
\lim_{n\gg k\gg m\rightarrow \infty }\mathbf{P}_{w}\left( \frac{1}{a_{m}}%
\left( S_{\tau _{r,n}}-S_{r}\right) \leq y|\mathcal{B}(ta_{k},n)\right) =%
\mathbf{P}\left( \min_{0\leq s\leq 1}Y_{s}\leq y\right) .
\end{equation*}
\end{theorem}

\textbf{Proof}. For fixed $w\geq 0,\,\varepsilon \in (0,1)$ and $N\in
\mathbb{N}$ we write
\begin{equation*}
\mathbf{P}_{w}\left( \frac{1}{a_{m}}\left( S_{\tau _{r,n}}-S_{n}\right) \leq
y;\mathcal{B}(ta_{k},n)\right) =J(r,k,n)+r_{m,n}(\varepsilon ,N),
\end{equation*}%
where
\begin{eqnarray*}
&&J(r,k,n) \\
&&\,=\mathbf{P}_{w}\left( \frac{1}{a_{m}}\left( S_{\tau _{r,n}}-S_{r}\right)
\leq y,\varepsilon a_{k}\leq S_{n}\leq ta_{k},\left\vert
S_{n}-S_{r}\right\vert \leq Na_{m},L_{n}\geq 0\right)
\end{eqnarray*}%
and
\begin{equation*}
0\leq r_{m,n}(\varepsilon ,N)\leq \mathbf{P}_{w}\left( 0\leq S_{n}\leq
\varepsilon a_{k},L_{n}\geq 0\right) +\mathbf{P}_{w}\left( \left\vert
S_{n}-S_{r}\right\vert >Na_{m},L_{n}\geq 0\right) .
\end{equation*}%
If $X_{1}\in \mathcal{D}(\alpha ,\beta )$, then, according to point 3) of
Lemma 8 in \cite{VDD2023} for any fixed $z\in (-\infty ,+\infty )$
\begin{equation*}
\mathbf{P}_{w}\mathbf{(}S_{n}-S_{r}\leq za_{m}|\mathcal{B}%
(ta_{k},n))\rightarrow \mathbf{P}\left( Y_{1}\leq z\right)
\end{equation*}%
as $n\gg k\gg m=n-r\rightarrow \infty $ uniformly in nonnegative $w$
satisfying the condition $w=o(a_{n})$. This statement combined with (\ref%
{End_point_neglogible}) shows that
\begin{equation*}
\lim_{\varepsilon \downarrow 0,N\uparrow \infty }\limsup_{n\gg k\gg
m\rightarrow \infty .}\frac{r_{m,n}(\varepsilon ,N)}{\mathbf{P}_{w}\mathbf{(}%
\mathcal{B}(ta_{k},n))}=0.
\end{equation*}%
Let us analyse the asymptotic behavior of the expression
\begin{eqnarray*}
&&J(m,k,n)=\int_{\varepsilon a_{k}-Na_{m}}^{ta_{k}+Na_{m}}\mathbf{P}%
_{w}\left( S_{r}\in dz,L_{r}\geq 0\right) \\
&&\qquad \times \mathbf{P}\left( \frac{1}{a_{m}}S_{\tau _{m}}\leq
y,\varepsilon a_{k}-z\leq S_{m}\leq ta_{k}-z,\left\vert S_{m}\right\vert
\leq Na_{m},S_{\tau _{m}}\geq -z\right)
\end{eqnarray*}%
as $n\gg k\gg m\rightarrow \infty $. First we note that, for sufficiently
large $k\gg m$ and $z\in \left( \varepsilon
a_{k}+Na_{m},ta_{k}-Na_{m}\right) $
\begin{equation*}
\left( -Na_{m},Na_{m}\right) \subset (\varepsilon a_{k}-z,ta_{k}-z).
\end{equation*}%
For $k\gg m$ meeting the restriction above we have
\begin{eqnarray*}
&&\mathbf{P}\left( \frac{1}{a_{m}}S_{\tau _{m}}\leq y,\varepsilon
a_{k}-z\leq S_{m}\leq ta_{k}-z,\left\vert S_{m}\right\vert \leq
Na_{m},S_{\tau _{m}}\geq -z\right) \\
&&\qquad \qquad =\mathbf{P}\left( \frac{1}{a_{m}}S_{\tau _{m}}\leq
y,\left\vert S_{m}\right\vert \leq Na_{m},S_{\tau _{m}}\geq -z\right) \\
&&\qquad \qquad =\mathbf{P}\left( \frac{1}{a_{m}}S_{\tau _{m}}\leq
y,\left\vert S_{m}\right\vert \leq Na_{m}\right) -h_{m,n}(z),
\end{eqnarray*}%
where the first term at the right-hand side of the last expression is
independent of $z$ and
\begin{equation*}
h_{m,n}(z)=\mathbf{P}\left( \frac{1}{a_{m}}S_{\tau _{m}}\leq y,\left\vert
S_{m}\right\vert \leq Na_{m},S_{\tau _{m}}<-z\right) .
\end{equation*}%
For random walk with $X_{1}\in \mathcal{D}(\alpha ,\beta )$ the
Donsker-Prokhorov invariance principle is valid. Therefore,
\begin{equation*}
\lim_{m\rightarrow \infty }\mathbf{P}\left( S_{\tau _{m}}\leq xa_{m}\right) =%
\mathbf{P}\left( \min_{0\leq s\leq 1}Y_{s}\leq x\right) ,x\in (-\infty
,+\infty ).
\end{equation*}%
Hence it follows that
\begin{eqnarray*}
\sup_{\varepsilon a_{k}-Na_{m}\leq z\leq ta_{k}+Na_{m}}h_{m,n}(z) &\leq
&\sup_{\varepsilon a_{k}-Na_{m}\leq z\leq ta_{k}+Na_{m}}\mathbf{P}\left(
S_{\tau _{m}}<-z\right) \\
&\leq &\mathbf{P}\left( \frac{S_{\tau _{m}}}{a_{m}}\leq N-\varepsilon \frac{%
a_{k}}{a_{m}}\right) \rightarrow 0
\end{eqnarray*}%
as $k\gg m\rightarrow \infty $. The estimates above allows us to conclude
that
\begin{eqnarray}
J(m,k,n) &=&\int_{\varepsilon a_{k}-Na_{m}}^{ta_{k}+Na_{m}}\mathbf{P}%
_{w}\left( S_{r}\in dz,L_{r}\geq 0\right)  \notag \\
&&\times \left( \mathbf{P}\left( \frac{1}{a_{m}}S_{\tau _{m}}\leq
y,\left\vert S_{m}\right\vert \leq Na_{m}\right) +o(1)\right)  \notag \\
&=&\mathbf{P}_{w}\left( S_{r}\in (\varepsilon
a_{k}+Na_{m},ta_{k}-Na_{m}),L_{r}\geq 0\right)  \notag \\
&&\times \mathbf{P}\left( \frac{1}{a_{m}}S_{\tau _{m}}\leq y,\left\vert
S_{m}\right\vert \leq Na_{m}\right)  \notag \\
&&+o\left( \mathbf{P}_{w}\left( S_{r}\leq ta_{k}+Na_{m},L_{r}\geq 0\right)
\right)  \notag
\end{eqnarray}%
as $n\gg k\gg m\rightarrow \infty $. Note that
\begin{eqnarray*}
\mathbf{P}_{w}\left( S_{r}\leq ta_{k}+Na_{m},L_{r}\geq 0\right) &\sim &%
\mathbf{P}_{w}\left( S_{n}\leq ta_{k}+Na_{m},L_{n}\geq 0\right) \\
&\sim &\mathbf{P}_{w}\left( \mathcal{B}(ta_{k},n)\right)
\end{eqnarray*}%
as $n\sim r\gg m$ in view of (\ref{BasicAsymptotic}), (\ref{AsympV}) and (%
\ref{Defb}). Further,
\begin{eqnarray*}
&&\mathbf{P}_{w}\left( S_{r}\in (\varepsilon
a_{k}+Na_{m},ta_{k}-Na_{m}),L_{r}\geq 0\right) \\
&&\qquad \qquad \qquad \times \mathbf{P}\left( \frac{1}{a_{m}}S_{\tau
_{m}}\leq y,\left\vert S_{m}\right\vert \leq Na_{m}\right) \\
&&\qquad \qquad =\mathbf{P}_{w}\left( S_{r}\leq ta_{k},L_{r}\geq 0\right)
\mathbf{P}\left( \frac{1}{a_{m}}S_{\tau _{m}}\leq y\right) +r_{m,n}^{\ast
}(\varepsilon ,N),
\end{eqnarray*}%
where
\begin{eqnarray*}
\left\vert r_{m,n}^{\ast }(\varepsilon ,N)\right\vert &\leq &\mathbf{P}%
_{w}\left( S_{r}\leq \varepsilon a_{k}+Na_{m},L_{r}\geq 0\right) \\
&&+\mathbf{P}_{w}\left( ta_{k}-Na_{m}\leq S_{r}\leq ta_{k},L_{r}\geq 0\right)
\\
&&+\mathbf{P}\left( \left\vert S_{m}\right\vert >Na_{m}\right) \mathbf{P}%
_{w}\left( S_{r}\leq ta_{k},L_{r}\geq 0\right) .
\end{eqnarray*}%
Using (\ref{IntrervalEstimate}), (\ref{AsympV}), (\ref{Regular1}), (\ref%
{Defb}) and (\ref{BasicAsymptotic}) we conclude that
\begin{eqnarray*}
\mathbf{P}_{w}\left( S_{r}\leq \varepsilon a_{k}+Na_{m},L_{r}\geq 0\right)
&\leq &Cb_{r}V^{-}(w)\int_{0}^{\varepsilon a_{k}+Na_{m}}V^{+}(u)du \\
&\leq &C_{1}b_{n}V^{-}(w)\varepsilon ^{\alpha \rho
+1}\int_{0}^{ta_{k}}V^{+}(u)du
\end{eqnarray*}%
for all sufficiently large $m$ and $n\gg k\gg m$. Thus,
\begin{equation}
\lim_{\varepsilon \downarrow 0,N\uparrow \infty }\limsup_{n\gg k\gg
m\rightarrow \infty }\frac{\mathbf{P}_{w}\left( S_{r}\leq \varepsilon
a_{k}+Na_{m},L_{r}\geq 0\right) }{\mathbf{P}_{w}\left( \mathcal{B}%
(ta_{k},n)\right) }=0.  \label{Close2}
\end{equation}%
It is not difficult to show by the same arguments that
\begin{eqnarray*}
&&\mathbf{P}_{w}\left( ta_{k}-Na_{m}\leq S_{r}\leq ta_{k},L_{r}\geq 0\right)
\\
&&\qquad \qquad \leq
Cb_{r}V^{-}(w)\int_{ta_{k}-Na_{m}}^{ta_{k}+Na_{m}}V^{+}(u)du \\
&&\qquad \qquad \leq C_{1}Nb_{n}V^{-}(w)a_{m}V^{+}(ta_{k}+Na_{m})=o\left(
b_{n}a_{k}V^{+}(ta_{k})\right) \\
&&\qquad \qquad =o\left( b_{n}V^{-}(w)\int_{0}^{ta_{k}}V^{+}(u)du\right)
\end{eqnarray*}%
and, therefore,
\begin{equation}
\limsup_{n\gg k\gg m\rightarrow \infty }\frac{\mathbf{P}_{w}\left(
ta_{k}-Na_{m}\leq S_{r}\leq ta_{k},L_{r}\geq 0\right) }{\mathbf{P}_{w}\left(
\mathcal{B}(ta_{k},n)\right) }=0.  \label{Close3}
\end{equation}

Since the distributions of the elements of the sequence $S_{m}/a_{m},%
\,m=1,2,...,$ converge, as $m\rightarrow \infty $ to the distribution of the
proper random variable $Y_{1}$ (see (\ref{dva})), it follows that
\begin{equation}
\lim_{N\rightarrow \infty }\mathbf{P}\left( \left\vert S_{m}\right\vert
>Na_{m}\right) =0.  \label{Close4}
\end{equation}%
Combining (\ref{Close2})--(\ref{Close4}) we conclude that
\begin{equation*}
\lim_{\varepsilon \downarrow 0,N\uparrow \infty }\limsup_{n\gg k\gg
m\rightarrow \infty .}\frac{\left\vert r_{m,n}^{\ast }(\varepsilon
,N)\right\vert }{\mathbf{P}_{w}\mathbf{(}\mathcal{B}(ta_{k},n))}=0.
\end{equation*}

This fact and the representation
\begin{eqnarray*}
&&\mathbf{P}_{w}\left( \frac{1}{a_{m}}\left( S_{\tau _{r,n}}-S_{r}\right)
\leq y;\mathcal{B}(ta_{k},n)\right) =J(m,k,n)+r_{m,n}(\varepsilon ,N) \\
&&\qquad \qquad \qquad \qquad =(1+o(1))\mathbf{P}_{w}\left( \mathcal{B}%
(ta_{k},n)\right) \mathbf{P}\left( \frac{1}{a_{m}}S_{\tau _{m}}\leq y\right)
\\
&&\qquad \qquad \qquad \qquad \qquad +r_{m,n}^{\ast }(\varepsilon
,N)+r_{m,n}(\varepsilon ,N)
\end{eqnarray*}%
imply that
\begin{eqnarray*}
\lim_{n\gg k\gg m\rightarrow \infty }\mathbf{P}_{w}\left( \frac{1}{a_{m}}%
\left( S_{\tau _{r,n}}-S_{r}\right) \leq y|\mathcal{B}(ta_{k},n)\right)
&=&\lim_{m\rightarrow \infty }\mathbf{P}\left( \frac{S_{\tau _{m}}}{a_{m}}%
\leq y\right) \\
&=&\mathbf{P}\left( \min_{0\leq s\leq 1}Y_{s}\leq y\right) .
\end{eqnarray*}

Theorem \ref{T_Minim} is proved.

\end{document}